\documentclass[11pt,a4paper]{article}
\usepackage[margin=1in]{geometry}

\usepackage{amsmath,amssymb,amsfonts}
\usepackage{amsthm}
\usepackage[title]{appendix}
\usepackage{xcolor}
\usepackage[hidelinks]{hyperref}

\theoremstyle{plain}
\newtheorem{theorem}{Theorem}[section]
\newtheorem{lemma}[theorem]{Lemma}
\newtheorem{corollary}[theorem]{Corollary}
\newtheorem{proposition}[theorem]{Proposition}
\theoremstyle{definition}
\newtheorem{definition}[theorem]{Definition}
\numberwithin{equation}{section}

\newcommand{\Id}{\mathrm{Id}}
\renewcommand{\Re}{\operatorname{Re}}
\renewcommand{\Im}{\operatorname{Im}}
\DeclareMathOperator{\supp}{supp}
\newcommand{\curl}{\mathrm{curl}}
\renewcommand{\div}{\mathrm{div}}
\newcommand{\comp}{0}
\newcommand{\compp}{\mathrm{comp}}
\newcommand{\loc}{\mathrm{loc}}

\begin{document}
\title{Spectral theory for transmission problems and resonances for transparent balls}
\author{Alden Waters\\
  {\small Leibniz Universit\"at Hannover}\\
  {\small Hannover, Germany}\\
  {\small \texttt{alden.waters@math.uni-hannover.de}}}
\date{}

\maketitle

\begin{abstract}
\normalsize
We construct generalized eigenfunctions and a spectral calculus for the transmission problem on Lipschitz domains. Then, focusing on the specific example of the ball, we prove a high-frequency resonance dichotomy for Maxwell's equations for a ball with positive, piecewise constant permittivity and permeability. If the product of these coefficients is no larger inside the ball than outside, the cutoff resolvent extends holomorphically to a strip of fixed width below the real axis, including the case of equal products. If the interior product is larger, we construct whispering-gallery resonances exponentially close to the real axis. The proofs use estimates for the transverse electric and transverse magnetic determinants that are uniform in angular momentum and complex frequency, together with a quasimode construction and Rouch\'e's theorem. Explicit generalized eigenfunctions identify these determinants with the scattering poles and give the spectral resolution and the normalization relating the scattering amplitude to the Mie T-matrix.
\end{abstract}



\section{Introduction and main theorems}

The main results of this paper are a resonance-free strip and the construction of exponentially narrow whispering-gallery resonances for Maxwell's equations for a transparent ball, although many of the spectral theory results are for Lipschitz domains. The permittivity $\epsilon$ and permeability $\mu$ are positive and piecewise constant, with subscripts $1$ and $0$ denoting the interior and exterior values. If $m_0=\sqrt{\epsilon_1\mu_1/(\epsilon_0\mu_0)}\leq1$, the cutoff resolvent is holomorphic in a strip of fixed width below the real axis at high frequency. If $m_0>1$, internal trapping produces resonances with exponentially small negative imaginary parts, so no such strip exists.

The main difficulty is to control the Maxwell transmission determinants uniformly in angular momentum and complex frequency. The explicit Mie formulas alone do not give this control: the interior and exterior Bessel arguments pass through different turning points, and both polarizations must be treated with their respective impedance factors. For $m_0\leq1$, we exclude zeros throughout the subturning, turning and oscillatory regimes, treating the endpoint $m_0=1$ separately. For $m_0>1$, the interior field can oscillate while the exterior field is evanescent; we construct real quasimodes in this regime and locate nearby resonances by Rouch\'e's theorem. This is the first time to the author's knowledge that the case of precise resonances-free strips and whispering-gallery modes for piecewise constant coefficient Maxwell's equations have been examined. The spectral theory and generalized eigenfunctions here are constructed for arbitrary Lipschitz domains, but for the main theorems we focus on the concrete example of the ball. 

For scalar transmission problems, the corresponding distinction between resonance-free regions and resonances approaching the real axis was established by Popov and Vodev \cite{PopovVodev,PopovVodevJE,PopovVodevAA} and Cardoso, Popov and Vodev \cite{CPV}. Related results include Galkowski's quantum Sabine law \cite{Galkowski2019Sabine} and the wavenumber-explicit estimates and resonance-free regions of Moiola and Spence \cite{MoiolaSpence}. For Maxwell's equations, the interface conditions also depend on the impedance, and the two polarizations lead to different scalar radial transmission problems. High-frequency Maxwell estimates in heterogeneous media are obtained in \cite{ChaumontFreletMoiolaSpence}; earlier asymptotic analyses of electromagnetic resonances in spheres and their perturbations include \cite{Schiller1993,GohsrichShahAiello2021,LamLeungYoung1992}.

To connect the determinant estimates to the resolvent, we construct normalized generalized eigenfunctions and derive the spectral resolution and resolvent jump. This extends the relation between scattering amplitudes and the T-matrix established for perfectly conducting spheres in \cite{FW}. Related spectral and scattering constructions appear in \cite{RBM,OS,NC,Parnovski,Galk}. The self-adjoint realization by a quadratic form provides the operator framework for these constructions.

Boundary integral treatments of Maxwell transmission include \cite{Torres1998,Rosen2019Spin,HelsingRosen2021Dirac}; in particular, Helsing and Ros{\'e}n use a Dirac formulation on Lipschitz domains. Costabel and Le Lou{\"e}r \cite{CostabelLeLouer2012ShapeII} study shape differentiability for dielectric scattering. For the scalar exterior Neumann problem, Galkowski, Marchand and Spence \cite{GalkowskiMarchandSpence2022Neumann} give frequency-explicit bounds for a second-kind boundary integral operator and its inverse.

For the operator framework and the generalized eigenfunctions, let $\Omega\subset\mathbb R^3$ be a bounded Lipschitz open set.  The isotropic material coefficients are
\begin{align}\label{eq:lipschitz-coefficients}
\epsilon(x)=\begin{cases}\epsilon_1,&x\in\Omega,\\\epsilon_0,&x\notin\overline\Omega,\end{cases}
\qquad
\mu(x)=\begin{cases}\mu_1,&x\in\Omega,\\\mu_0,&x\notin\overline\Omega.\end{cases}
\end{align}
The values on $\partial\Omega$ do not affect the operators. The two main resonance theorems and the explicit transmission coefficients below specialize to a ball.

We also assume that $0<\epsilon_j<\infty$ and $0<\mu_j<\infty$, $j=0,1$. We then consider the spectral version of Maxwell's equations for the electric field $\mathcal{E}$ and the magnetic field $\mathcal{H}$ as:
\begin{align}\label{spectralMax}
&\curl \mathcal{E}-\mathrm{i}\lambda \mu \mathcal{H}=0\quad \mathrm{in} \quad \mathbb{R}^3 \\& \nonumber
\div (\epsilon\cdot \mathcal{E})=0\\& \nonumber
\curl \mathcal{H}+\mathrm{i}\lambda \epsilon \mathcal{E}=0\\& \nonumber
\div (\mu\cdot \mathcal{H})=0\\& \nonumber
\nu\times \mathcal{E}_{\mathrm{int}}= \nu\times \mathcal{E}_{\mathrm{ext}} \quad\mathrm{on}\quad \partial\Omega\\& \nonumber
\nu\times \mathcal{H}_{\mathrm{int}}=\nu \times \mathcal{H}_{\mathrm{ext}} \quad \mathrm{on} \quad \partial\Omega. \nonumber
\end{align}
The Hilbert space for the electric field is
\begin{align}
(L^2(\mathbb{R}^3;\epsilon(x)\,\mathrm{d}x))^3 \quad\quad \langle u,v\rangle_{(L^2(\mathbb{R}^3;\epsilon(x)\,\mathrm{d}x))^3}=\int\limits_{\mathbb{R}^3}\epsilon u\cdot \overline{v}\,\mathrm{d}x.
\end{align}
We abbreviate this as $(L^2(\epsilon\,\mathrm{d}x))^3$.
Restricting to the divergence-free subspace we have
\begin{align}
L_{\epsilon}(\div \,0)=\{u\in (L^2(\epsilon\,\mathrm{d}x))^3: \nabla\cdot(\epsilon u)=0\}.
\end{align}
It follows that the operator $P_E$ associated to the Maxwell system electric field acting on 3 dimensional vectors $u$ is
\begin{align}\label{system}
P_Eu=\epsilon^{-1}\nabla\times \mu^{-1}\nabla\times u.
\end{align}
The boundary conditions for the electric field for this system can also be recast as:
\begin{align}\label{transmission}
&\nu\times \mathcal{E}_{\mathrm{int}}=\nu\times \mathcal{E}_{\mathrm{ext}} \\& \nonumber
\nu\times \mu_1^{-1}\nabla\times \mathcal{E}_{\mathrm{int}}=\nu\times \mu_0^{-1}\nabla\times \mathcal{E}_{\mathrm{ext}} \quad \text{on }\partial\Omega.
\end{align}
We use $H(\curl,D)=\{u\in L^2(D;\mathbb C^3):\curl u\in L^2(D;\mathbb C^3)\}$, $D=\Omega,\Omega^c$ or $\mathbb{R}^3$, and the weighted orthogonal Helmholtz decomposition
\begin{align}\label{HD}
(L^2(\epsilon\,\mathrm{d}x))^3
=L_\epsilon(\div0)\oplus\nabla\dot H^1(\mathbb R^3).
\end{align}
Here $\dot H^1(\mathbb R^3)$ is the homogeneous Sobolev space with square-integrable gradient. The transmission conditions are encoded by the domain
\begin{align*}
\mathcal D(P_E)=\{u\in L_\epsilon(\div0)\cap H(\curl,\mathbb R^3):
\mu^{-1}\curl u\in H(\curl,\mathbb R^3)\}.
\end{align*}
Indeed, global $H(\curl)$ regularity is equivalent to matching tangential traces across the interface; see \cite{KH}. Let $\nu$ be the outward pointing normal vector to $\Omega$. We define the jump as $
[[F]]=(F^{\mathrm{int}}-F^{\mathrm{ext}})|_{\partial\Omega}.$
For a piecewise smooth field $G$, with components $G_0,G_1$, the curl Green identity with a compactly supported test field $v$ reads
\begin{align}\label{byparts}
&\int_{\Omega}\curl G_1\cdot\overline v\,\mathrm{d}x
+\int_{\mathbb R^3\setminus\overline\Omega}\curl G_0\cdot\overline v\,\mathrm{d}x=\int_{\mathbb R^3}G\cdot\overline{\curl v}\,\mathrm{d}x
+\langle[[\nu\times G]],v_T\rangle.
\end{align}
Here $v_T$ is the tangential trace; the identity extends to piecewise $H(\curl)$ fields by the $H^{-1/2}(\operatorname{Div}_{\partial\Omega})$--$H^{-1/2}(\operatorname{Curl}_{\partial\Omega})$ trace duality, defined by this Green identity, \cite{KH}. All boundary pairings below use this duality, conjugate-linear in the second argument.
\begin{proposition}
The operator $P_E$ is a non-negative self-adjoint operator on $L_{\epsilon}(\div 0)$, i.e., let $V_E$ be the space
  \begin{align}
V_E=\{u\in (L^2(\mathbb{R}^3;\epsilon(x)\,\mathrm{d}x))^3: \curl u\in (L^2(\mathbb{R}^3))^3 \quad \mathrm{div}(\epsilon u)=0\}
\end{align}
and we define the corresponding quadratic form as:
\begin{align}
q_E(u,v)=\int\limits_{\mathbb{R}^3}\mu^{-1}(\curl u)\cdot \overline{(\curl v)}\,\mathrm{d}x.
\end{align}
Then $q_E$ is a densely defined, closed, non-negative symmetric form on $L_{\epsilon}(\mathrm{div}\, 0)$ and the self adjoint operator $P_E$, is realized by $q_E(u,u)=|| P^{\frac{1}{2}}_Eu ||^2_{(L^2(\epsilon\,\mathrm{d}x))^3}$.
\end{proposition}
\begin{proof}
The form is symmetric and non-negative. The squared form norm $q_E(u,u)+\|u\|_{L^2(\epsilon\,\mathrm{d}x)^3}^2$ is equivalent to the squared $H(\curl)$ norm, and the divergence constraint is closed, so $q_E$ is closed. Let $\Pi$ be the weighted orthogonal projection in \eqref{HD}. Since $\Pi v=v-\nabla p$ and $\curl\Pi v=\curl v$, it maps $C_0^\infty$ into $V_E$. Boundedness of $\Pi$ and density of $C_0^\infty$ in $L^2$ show that $V_E$ is dense in $L_\epsilon(\div0)$.

The representation theorem \cite[Theorem VIII.15]{RS} gives a non-negative self-adjoint operator $A$. For $u\in\mathcal D(P_E)$, the global curl Green identity above gives $q_E(u,v)=\langle P_Eu,v\rangle_{L^2(\epsilon\,\mathrm{d}x)^3}$ for all $v\in V_E$, hence $P_E\subset A$. Conversely, if $u\in\mathcal D(A)$ and $g=Au$, then
\begin{equation*}
q_E(u,v)=\langle g,v\rangle_{L^2(\epsilon\,\mathrm{d}x)^3},\qquad v\in V_E.
\end{equation*}
Both sides vanish on gradients, so \eqref{HD} extends this identity to all $v\in H(\curl,\mathbb R^3)$. Thus $\curl(\mu^{-1}\curl u)=\epsilon g$ distributionally on $\mathbb R^3$. Since $\mu^{-1}\curl u\in L^2$, it belongs to $H(\curl)$, so $u\in\mathcal D(P_E)$ and $P_Eu=g$. Therefore $A=P_E$.
\end{proof}
\begin{corollary}
The magnetic operator $P_{H}=\mu^{-1}\curl \epsilon^{-1}\curl$ on $L_{\mu}(\div 0)$ with corresponding quadratic form $q_H(u,v)=\int\limits_{\mathbb{R}^3}\epsilon^{-1}\curl u\cdot\overline{\curl v}\,\mathrm{d}x$ is self adjoint.
\end{corollary}
\begin{proof}
Interchange $\epsilon$ and $\mu$ in the preceding proof. The form domain is $V_H=H(\curl,\mathbb R^3)\cap L_\mu(\div0)$.
\end{proof}
Let $\Pi_E$ and $\Pi_H$ be the weighted orthogonal projections onto $L_\epsilon(\div0)$ and $L_\mu(\div0)$. Non-negativity and self-adjointness give
\begin{align}
R_{\lambda,E}=(P_E-\lambda^2)^{-1},\qquad
R_{\lambda,H}=(P_H-\lambda^2)^{-1},\qquad \Im\lambda>0.
\end{align}
For real $\lambda\ne0$, we use their boundary values from the upper half-plane, established below, as operators from compactly supported data to $L^2_{\loc}$. To allow scalar cutoffs, we extend the electric resolvent to all of $(L^2(\epsilon\,\mathrm{d}x))^3$ by
\begin{align}
\tilde R_{\lambda,E}=R_{\lambda,E}\Pi_E-\lambda^{-2}(\Id-\Pi_E).
\end{align}
The second term is the resolvent on the gradient summand, where the curl--curl operator vanishes. The magnetic extension is defined analogously. We make the following definition which is adapted from \cite{DZbook} for our setting: 
\begin{definition}\label{resonanceP}
A vector $u\in H_{loc}(\curl)$ is a \textbf{resonance state} corresponding to $\lambda\in\mathbb{C}\setminus\{0\}$ if and only if $(P_E-\lambda^2)u=0$ (resp. $(P_H-\lambda^2)u=0$) and there exists $g\in (L^2_{\comp}(\mathbb{R}^3))^3$, $R>0$ with $\overline\Omega\subset B_R$ such that
\begin{align}
u|_{\mathbb{R}^3\setminus B(0,R)}=R_{0,\lambda}g|_{\mathbb{R}^3\setminus B(0,R)}
\end{align}
i.e. $u$ is \textit{purely outgoing}. Here $R_{0,\lambda}=(P_0-\lambda^2)^{-1}$ denotes the free constant coefficient $\epsilon_1=\epsilon_0$, $\mu_1=\mu_0$ outgoing resolvent.
\end{definition}
For the following two theorems we take $\Omega=B_{r_0}=\{x\in\mathbb R^3:|x|<r_0\}$ and $m_0=\sqrt{\epsilon_1\mu_1/(\epsilon_0\mu_0)}$. The two cases for $m_0$ give the following dichotomy.
\begin{theorem}\label{maintheorem}
Let $0<m_0\leq 1$ then there exist positive $C_1,C_2$ depending on $\mu_0,\mu_1,\epsilon_0,\epsilon_1$ and $r_0$ such that for every smooth compactly supported (scalar) function $\chi$ which is one in a neighborhood of the ball, $\chi \tilde{R}_{\lambda,E}\chi$ and $\chi \tilde{R}_{\lambda,H}\chi$ are holomorphic in
\begin{align}
\left\{ \lambda\in \mathbb{C}: \,\, |\mathrm{Re}(\lambda)|\geq C_1, \quad \Im \lambda\geq -C_2 \right\}.
\end{align}
\end{theorem}
Furthermore in the case in which $\epsilon_1\mu_1>\epsilon_0\mu_0$ we have resonances, which manifest themselves in the following theorem.
\begin{theorem}\label{maintheorem2}
Let $1<m_0$, then no resonance free strip exists. There exist constants $C,c>0$ which depend on the material parameters $m_0,\mu_0,\mu_1,\epsilon_0,\epsilon_1$ and $r_0$ as well as a sequence of resonances $\{\lambda_j\}_{j\in\mathbb{N}}$ of the operator $P_E$ such that
\begin{align}
\mathrm{Re}\lambda_j\rightarrow\infty \quad \mathrm{Im}\lambda_j<0 \quad |\mathrm{Im}\lambda_j|\leq C\mathrm{e}^{-c\mathrm{Re}\lambda_j}.
\end{align}
\end{theorem}
\section{Vector spherical harmonics and the T-matrix}\label{sec:tmatrix}
Let $Y_{\ell m}$, $\ell\geq0$, $|m|\leq\ell$, be a real orthonormal basis of scalar spherical harmonics, extended homogeneously of degree zero to $x=r\hat r$. For $\ell\geq1$ set
\begin{gather}
\Psi_{1,\ell,m}=\frac{\nabla Y_{\ell m}\times x}{\sqrt{\ell(\ell+1)}},\qquad
\Psi_{2,\ell,m}=\frac{r\nabla Y_{\ell m}}{\sqrt{\ell(\ell+1)}},\qquad
\Psi_{3,\ell,m}=\hat rY_{\ell m}.
\end{gather}
The first two families form an orthonormal basis of
\begin{align}
L^2_t(\mathbb S^2;\mathbb C^3)=\{F\in L^2(\mathbb S^2;\mathbb C^3):\hat r\cdot F=0\}.
\end{align}
Together with the radial family, including $\Psi_{3,0,0}$, they form an orthonormal basis of $L^2(\mathbb S^2;\mathbb C^3)$. We write $Y_n=\Psi_{j,\ell,m}$ with $n=(j,\ell,m)$ on the unit sphere. The components of $\Psi_{1,\ell,m}$ have harmonic degree $\ell$, whereas those of $\Psi_{2,\ell,m}$ and $\Psi_{3,\ell,m}$ combine degrees $\ell-1$ and $\ell+1$; see \cite[Theorem 2.46 and Corollary 2.47]{KH}.

For $z_\ell=j_\ell,h_\ell^{(1)}$ or $h_\ell^{(2)}$, the Hansen fields are
\begin{align*}
M_{\ell m}(kr)&=\curl(xz_\ell(kr)Y_{\ell m})
=\sqrt{\ell(\ell+1)}z_\ell(kr)\Psi_{1,\ell,m},\\
N_{\ell m}(kr)&=k^{-1}\curl M_{\ell m}(kr)\\
&=\frac{\ell(\ell+1)}{kr}z_\ell(kr)\Psi_{3,\ell,m}
+\frac{\sqrt{\ell(\ell+1)}}{kr}(yz_\ell(y))'\big|_{y=kr}\Psi_{2,\ell,m}.
\end{align*}
They satisfy $\curl M=kN$, $\curl N=kM$ and $\div M=\div N=0$. Set
\begin{align}
k_j=\lambda\sqrt{\epsilon_j\mu_j},\qquad \eta_j=\sqrt{\frac{\mu_j}{\epsilon_j}},\qquad j=0,1.
\end{align}
Either field with $k=k_j$ gives an electric solution in medium $j$, with magnetic field $H=(\mathrm{i}\lambda\mu_j)^{-1}\curl E$. These identities follow from the scalar Bessel equation; see \cite[Lemma 2.42 and Theorem 2.43]{KH} and \cite{scat}. Unlike the PEC construction in \cite{FW}, the transmission eigenfunctions are defined on all of $\mathbb R^3$.

We now define the generalized eigenfunctions piecewise for $B_{r_0}$. Here and below $x_j=k_jr_0$ for $j=0,1$.  The superscript $0$ refers to $j_{\ell}$ the superscript $-$ refers to $h_{\ell}^{(1)}$ in the following:
\begin{align}\label{piecewiseE}
E_{\lambda,1}(Y_n)=k_0\sqrt{\frac{2\eta_0}{\pi\ell(\ell+1)}}\begin{cases} t_{\ell}^{1}(\lambda) M_{\ell m}^{0}(k_1 r) \quad r<r_0\\
M_{\ell m}^{0}(k_0 r)+A_{\ell}^{1}(\lambda)M_{\ell m}^{-}(k_0 r) \quad r>r_0
\end{cases}
\end{align}
and
\begin{align}\label{piecewiseH}
E_{\lambda,2}(Y_n)=k_0\sqrt{\frac{2\eta_0}{\pi\ell(\ell+1)}}\begin{cases} t_{\ell}^{2}(\lambda) N_{\ell m}^{0}(k_1 r) \quad r<r_0\\
N_{\ell m}^{0}(k_0 r)+A_{\ell}^{2}(\lambda)N_{\ell m}^{-}(k_0 r) \quad r>r_0.
\end{cases}
\end{align}
Here $n=(j,\ell,m)$ in $E_{\lambda,j}(Y_n)$. The prefactors give delta normalization in $L^2(\epsilon\,\mathrm{d}x)^3$, with orthogonal polarizations; this is established in Section~\ref{sec:geneig}.
Set $\nu=\ell+\frac12$ and introduce the Riccati--Bessel functions
\begin{align}\label{eq:riccati-definitions}
\psi_{\ell}(y)&=yj_{\ell}(y)=\sqrt{\frac{\pi y}{2}}J_\nu(y),\nonumber\\
\xi_{\ell}(y)&=yh_{\ell}^{(1)}(y)=\sqrt{\frac{\pi y}{2}}H_\nu^{(1)}(y).
\end{align}
Primes denote differentiation with respect to the full argument.
\subsection{The transverse electric coefficient \texorpdfstring{$A_{\ell}^1$}{A1}}
Tangential continuity of $\mathcal E$ and $\mathcal H$ gives the TE matching system \cite{scat}
\begin{align}
&t_{\ell}^1j_{\ell}(x_1)=j_{\ell}(x_0)+A_{\ell}^1h_{\ell}^{(1)}(x_0)\\&
\frac{t_{\ell}^1}{\mu_1}\psi'_{\ell}(x_1)=\frac{1}{\mu_0}\left(\psi_{\ell}'(x_0)+A_{\ell}^1\xi_{\ell}'(x_0)\right).
\end{align}
It follows that
\begin{align}
A_{\ell}^1(\lambda)&=\frac{\eta_1\psi_{\ell}(x_1)\psi'_{\ell}(x_0)-\eta_0\psi_{\ell}(x_0)\psi_{\ell}'(x_1)}{\eta_0\xi_{\ell}(x_0)\psi_{\ell}'(x_1)-\eta_1\psi_{\ell}(x_1)\xi_{\ell}'(x_0)}.
\end{align}
The interior amplitude is
\begin{align}
t_{\ell}^1(\lambda)=-\frac{\mathrm{i}\mu_1}{x_0}\frac{1}{\mu_0h_{\ell}^{(1)}(x_0)\psi_{\ell}'(x_1)-\mu_1j_{\ell}(x_1)\xi_{\ell}'(x_0)}.
\end{align}
The interior amplitude follows from the Wronskian \eqref{wronskian}.

\subsection{The transverse magnetic coefficient \texorpdfstring{$A_{\ell}^2$}{A2}}
The TM matching system is
\begin{align}
&\frac{t_{\ell}^2}{k_1}\psi'_{\ell}(x_1)=\frac{1}{k_0}\left(\psi_{\ell}'(x_0)+A_{\ell}^2\xi_{\ell}'(x_0)\right) \\&
\frac{k_1 t_{\ell}^2}{\mu_1}j_{\ell}(x_1)=\frac{k_0}{\mu_0}\left(j_{\ell}(x_0)+A_{\ell}^2h_{\ell}^{(1)}(x_0)\right).
\end{align}
We solve the system to obtain
\begin{align}
A_{\ell}^2(\lambda)&=\frac{\eta_1\psi_{\ell}(x_0)\psi'_{\ell}(x_1)-\eta_0\psi_{\ell}(x_1)\psi_{\ell}'(x_0)}{\eta_0\psi_{\ell}(x_1)\xi_{\ell}'(x_0)-\eta_1\xi_{\ell}(x_0)\psi_{\ell}'(x_1)}
\end{align}
and also
\begin{align}
t_{\ell}^2(\lambda)=\frac{\frac{\mathrm{i}k_0k_1\mu_1}{x_0}}{k_1^2\mu_0j_{\ell}(x_1)\xi'_{\ell}(x_0)-k_0^2\mu_1h_{\ell}^{(1)}(x_0)\psi_{\ell}'(x_1)}.
\end{align}
\begin{lemma}
For $\lambda>0$, let $D_\ell^j(\lambda)\ne0$ be the denominators of $A_\ell^j(\lambda)$. Define $S(\lambda)$ diagonally in the vector spherical harmonic basis by
\begin{align}
S(\lambda)Y_{j,\ell,m}=S_{\ell}^j(\lambda)Y_{j,\ell,m} \quad\quad S_{\ell}^j(\lambda)=1+2A_{\ell}^j(\lambda)=-\frac{\overline{D_{\ell}^j(\lambda)}}{D_{\ell}^j(\lambda)}.
\end{align}
In particular, $S(\lambda)$ is unitary.
\end{lemma}
\begin{proof}
Write $A_\ell^j=N_\ell^j/D_\ell^j$ using the preceding formulas. For real $\lambda>0$, the numerators $N_\ell^j$ are real and $\xi_\ell+\overline{\xi_\ell}=2\psi_\ell$, with the same identity for their derivatives. Hence $D_\ell^j+\overline{D_\ell^j}=-2N_\ell^j$, which gives the stated identity and unitarity.
\end{proof}
We set
\begin{align}
\tilde{j}_{\lambda}(Y_n)(k_ 0r,\theta,\phi)=k_0\sqrt{\frac{2\eta_0}{\pi\ell(\ell+1)}}M_{\ell m}^{0}(k_0 r)
\end{align}
and
\begin{align}
\tilde{H}_{\lambda}(Y_n)=\frac{k_0}{\mathrm{i}\eta_0}\sqrt{\frac{2\eta_0}{\pi\ell(\ell+1)}}M_{\ell m}^0(k_0r).
\end{align}
These free fields are tangential to spheres where $\Omega=B_{r_0}$ and satisfy the weighted divergence constraints. We now give the construction for a general bounded Lipschitz $\Omega$, retaining $R_{\lambda,E}$ on $L_\epsilon(\div0)$ and $R_{\lambda,H}$ on $L_\mu(\div0)$. We encode the coefficient contrast by compactly supported functionals on the form domains.
\begin{definition}
We define the compactly supported functional $\Theta_{\lambda,E}(Y_n)$ on the space $V_E$ by
\begin{align}
\Theta_{\lambda, E}(Y_n)[\varphi]&=\lambda^2\int\limits_{\Omega}\left( \frac{\epsilon_0\mu_0}{\mu_1}-\epsilon_1\right)\tilde{j}_{\lambda}(Y_n)\cdot \overline{\varphi}\,\mathrm{d}x-\langle [[\nu\times \mu^{-1}\curl\tilde{j}_{\lambda}(Y_n)]],\varphi_T\rangle_{\partial\Omega}
\end{align}
and also the functional $\Theta_{\lambda,H}(Y_n)$ on $V_H$ by
\begin{align}
\Theta_{\lambda,H}(Y_n)[\varphi]&=\lambda^2\int\limits_{\Omega}\left( \frac{\epsilon_0\mu_0}{\epsilon_1}-\mu_1\right)\tilde{H}_{\lambda}(Y_n)\cdot \overline{\varphi}\,\mathrm{d}x-\langle [[\nu\times \epsilon^{-1}\curl \tilde{H}_{\lambda}(Y_n)]],\varphi_T\rangle_{\partial\Omega}.
\end{align}
Give $V_E$ the form norm $\|u\|_{V_E}^2=q_E(u,u)+\|u\|_{L^2(\epsilon \mathrm{d}x)^3}^2$ and write $V_E^*$ for its continuous anti-dual. Then for $\Im\lambda>0$,
$R_{\lambda,E}:V_E^*\longrightarrow V_E$.
For $\Im\lambda>0$, $u=R_{\lambda,E}\Theta_{\lambda,E}(Y_n)\in V_E$ is the unique solution of
\begin{equation*}
q_E(u,\varphi)-\lambda^2\langle u,\varphi\rangle_{(L^2(\epsilon \mathrm{d}x))^3}=\Theta_{\lambda,E}(Y_n)[\varphi]\qquad\forall\varphi\in V_E,
\end{equation*}
that is, $\Theta_{\lambda,E}(Y_n)=(P_E-\lambda^2)u$ in the weak sense. Equivalently,
\begin{equation*}
\langle R_{\lambda,E}\Theta_{\lambda,E}(Y_n),f\rangle_{(L^2(\epsilon \mathrm{d}x))^3}=\Theta_{\lambda,E}(Y_n)\big[(P_E-\bar\lambda^2)^{-1}f\big],\qquad f\in L_\epsilon(\div 0).
\end{equation*}
For $\lambda>0$ we define $R_{\lambda,E}\Theta_{\lambda,E}(Y_n)$ by
\begin{equation*}
\langle R_{\lambda,E}\Theta_{\lambda,E}(Y_n),f\rangle_{(L^2(\epsilon \mathrm{d}x))^3}=\Theta_{\lambda,E}(Y_n)\big[R_{-\lambda,E}f\big],\qquad f\in L^2_{\compp}\cap L_\epsilon(\div 0);
\end{equation*}
Here $R_{\lambda,E}\Theta_{\lambda,E}$ denotes the $H_{\loc}(\curl)$ limit from $\Im\lambda>0$, justified in the proof below; the displayed pairing is its duality characterization. 
\end{definition}
For the general interface define $p_{\lambda,E}\in\dot H^1(\mathbb R^3)$ by
\begin{equation}\label{eq:incident-gradient}
\int_{\mathbb R^3}\epsilon\nabla p_{\lambda,E}\cdot\overline{\nabla v}\,\mathrm{d}x
=(\epsilon_1-\epsilon_0)\int_\Omega\tilde j_\lambda(Y_n)\cdot\overline{\nabla v}\,\mathrm{d}x,
\qquad v\in\dot H^1(\mathbb R^3).
\end{equation}
Uniform positivity gives a unique solution by Lax--Milgram. Define $p_{\lambda,H}$ by interchanging $\epsilon,\mu$ and replacing $\tilde j_\lambda$ by $\tilde H_\lambda$. Then $\div(\epsilon(\tilde j_\lambda-\nabla p_{\lambda,E}))=0$, and similarly for the magnetic field. For the ball the right-hand sides vanish by the fact that  $\div(\tilde{j}_{\lambda})=0$ and $\tilde{j}_{\lambda}\cdot \nu=0$, so correspondingly $p_{\lambda,E}=p_{\lambda,H}=0$.
\begin{lemma}\label{resolventdist}
For compactly supported test fields $\varphi\in H(\curl,\mathbb{R}^3)$
\begin{align*}
\Theta_{\lambda,E}[\varphi]=q_E(\tilde{j}_{\lambda},\varphi)-\lambda^2\langle \tilde{j}_{\lambda},\varphi\rangle_{(L^2(\epsilon \mathrm{d}x))^3}.
\end{align*}
For general $\varphi\in V_E$, the left-hand side is defined by the preceding volume-and-boundary formula; the magnetic identity below is interpreted in the same way on $V_H$. For general $\Omega$ we define, for $\Im\lambda>0$ and by outgoing limits for $\lambda>0$,
\begin{align*}
E_{\lambda,1}(Y_n)=\tilde{j}_{\lambda}(Y_n)-\nabla p_{\lambda,E}-R_{\lambda,E}\Theta_{\lambda,E}(Y_n).
\end{align*}
This field solves the transmission problem with incident field $\tilde j_\lambda(Y_n)$ and outgoing scattered part $-\nabla p_{\lambda,E}-R_{\lambda,E}\Theta_{\lambda,E}(Y_n)$. For $\Omega=B_{r_0}$ it agrees with \eqref{piecewiseE}; explicitly,
\begin{align}\label{casesR}
R_{\lambda,E}\Theta_{\lambda,E}(Y_n)=\begin{cases}k_0\sqrt{\frac{2\eta_0}{\pi\ell(\ell+1)}}(M_{\ell m}^0(k_0r)-t_{\ell}^{1}(\lambda)M_{\ell m}^0(k_1 r)) \quad r<r_0\\
-k_0\sqrt{\frac{2\eta_0}{\pi\ell(\ell+1)}}A_{\ell}^1(\lambda)M_{\ell m}^-(k_0r) \quad r>r_0.\end{cases}
\end{align}
Similarly, for general $\Omega$ define
\begin{align*}
H_{\lambda}(Y_n)=\tilde{H}_{\lambda}(Y_n)-\nabla p_{\lambda,H}-R_{\lambda,H}\Theta_{\lambda,H}(Y_n)
\end{align*}
where
\begin{align}
\Theta_{\lambda,H}[\varphi]=q_{H}(\tilde{H}_{\lambda},\varphi)-\lambda^2\langle\tilde{H}_{\lambda},\varphi\rangle_{(L^2(\mu \mathrm{d}x))^3}
\end{align}
and define the second electric family by
\begin{align}
E_{\lambda,2}(Y_n)=\frac{\mathrm{i}}{\lambda}\epsilon^{-1}\curl H_{\lambda}(Y_n).
\end{align}
These fields satisfy the corresponding transmission problems and agree with the ball construction, in particular \eqref{piecewiseH}, when $\Omega=B_{r_0}$.
\end{lemma}
\begin{proof}
Fix an incident mode $Y_n$ and initially let $\Im\lambda>0$. We first verify the electric construction for a general bounded Lipschitz inclusion. Integration by parts shows that $\Theta_{\lambda,E}$ is the residual of the free incident field $\tilde j_\lambda$ in the transmission equation. More explicitly, as a functional on $H(\curl,\mathbb R^3)$ it is
\begin{equation*}
\Theta_{\lambda,E}[v]=\int_{\mathbb R^3}F_\lambda\cdot\overline v\,\mathrm{d}x
+\int_{\mathbb R^3}G_\lambda\cdot\overline{\curl v}\,\mathrm{d}x,
\end{equation*}
where
\begin{equation*}
F_\lambda=-\lambda^2(\epsilon_1-\epsilon_0)\mathbf1_\Omega\tilde j_\lambda,
\qquad
G_\lambda=(\mu_1^{-1}-\mu_0^{-1})\mathbf1_\Omega\curl\tilde j_\lambda.
\end{equation*}
This agrees with the volume-and-boundary definition and gives the stated weak identity for compactly supported test fields.

Set $u=R_{\lambda,E}\Theta_{\lambda,E}\in V_E$. The resolvent equation is initially tested only on $V_E$, so we must account for the gradient part of an arbitrary test field. Since $(\Id-\Pi_E)v$ is a gradient, the definition of $p_{\lambda,E}$ in \eqref{eq:incident-gradient} gives
\begin{align*}
\Theta_{\lambda,E}[(\Id-\Pi_E)v]
&=-\lambda^2\langle\nabla p_{\lambda,E},v\rangle_{L^2(\epsilon\,\mathrm{d}x)^3}.
\end{align*}
Testing the resolvent equation on $\Pi_Ev\in V_E$ therefore yields, for every $v\in H(\curl,\mathbb R^3)$,
\begin{align*}
q_E(u,v)-\lambda^2\langle u,v\rangle_{L^2(\epsilon\,\mathrm{d}x)^3}
&=\Theta_{\lambda,E}[\Pi_Ev]
=\Theta_{\lambda,E}[v]+\lambda^2\langle\nabla p_{\lambda,E},v\rangle_{L^2(\epsilon\,\mathrm{d}x)^3}.
\end{align*}
Equivalently, on all of $\mathbb R^3$ in the distributional sense,
\begin{equation*}
\curl\mu^{-1}\curl u-\lambda^2\epsilon u
=F_\lambda+\curl G_\lambda+\lambda^2\epsilon\nabla p_{\lambda,E}.
\end{equation*}
The residual of $\tilde j_\lambda$ is $F_\lambda+\curl G_\lambda$, whereas that of $\nabla p_{\lambda,E}$ is $-\lambda^2\epsilon\nabla p_{\lambda,E}$. Subtracting these identities proves that
\begin{equation*}
E:=\tilde j_\lambda-\nabla p_{\lambda,E}-u
\quad\hbox{satisfies}\quad
\curl\mu^{-1}\curl E-\lambda^2\epsilon E=0.
\end{equation*}
Moreover, \eqref{eq:incident-gradient} and $u\in V_E$ give $\div(\epsilon E)=0$. We have $E\in H_{\loc}(\curl)$, and its equation gives $\mu^{-1}\curl E\in H_{\loc}(\curl)$. These two global regularity statements imply the two tangential transmission conditions. The scattered field $E-\tilde j_\lambda=-\nabla p_{\lambda,E}-u$ solves the free homogeneous equation outside $\overline\Omega$ and is square integrable there for $\Im\lambda>0$, hence is outgoing.

It remains to justify the boundary value at real $\lambda>0$. Introduce the auxiliary field
\begin{equation*}
H_u:=\frac{1}{\mathrm{i}\lambda}\big(\mu^{-1}\curl u-G_\lambda\big).
\end{equation*}
The equation for $u$ is then equivalent to the first-order system
\begin{align*}
\curl u-\mathrm{i}\lambda\mu H_u&=\mu G_\lambda,\\
\curl H_u+\mathrm{i}\lambda\epsilon u
&=\frac{F_\lambda}{\mathrm{i}\lambda}-\mathrm{i}\lambda\epsilon\nabla p_{\lambda,E}.
\end{align*}
The fields $F_\lambda,G_\lambda$ are compactly supported and square integrable. The potential $p_{\lambda,E}$ is harmonic outside a ball and has zero total flux, because the source in its divergence-form equation is the divergence of a compactly supported vector field. Its finite-energy exterior harmonic expansion therefore gives $\nabla p_{\lambda,E}=O(r^{-3})$. Both sources in the first-order system consequently belong to $L^2_s=\{f:(1+|x|^2)^{s/2}f\in L^2\}$ for $1/2<s<3/2$.

We may now apply the whole-space limiting absorption principle \cite[Theorem 3.10 and Corollary 3.12]{OsterbrinkPauly2019LAP}, also summarized in \cite[Theorem 4.3]{OsterbrinkPauly2020}. The coefficients are uniformly positive and bounded and are constant outside a ball; local Maxwell compactness follows by localization to balls. Theorem~\ref{Rellich} rules out the homogeneous radiating kernel at every real $\lambda>0$. We obtain the outgoing $H_{\loc}(\curl)$ limit of $u$, with bounds locally uniform on compact subintervals of $(0,\infty)$. This also gives the outgoing boundary value of $E-\tilde j_\lambda$. For compactly supported volume data $f$, the gradient part of its weighted Helmholtz decomposition likewise decays as $O(r^{-3})$, so the same argument applies to $R_{\lambda,E}\Pi_Ef$.

Interchanging $\epsilon$ and $\mu$ and using the incident field $\tilde H_\lambda$ gives the total magnetic field
\begin{equation*}
H_\lambda=\tilde H_\lambda-\nabla p_{\lambda,H}-R_{\lambda,H}\Theta_{\lambda,H}.
\end{equation*}
This is the field used to construct the second electric family; $H_u$ above was only an auxiliary field for the limiting absorption argument. Since $\curl\epsilon^{-1}\curl H_\lambda=\lambda^2\mu H_\lambda$, setting $E_{\lambda,2}=\mathrm{i}\epsilon^{-1}\curl H_\lambda/\lambda$ gives both Maxwell curl equations and $\div(\epsilon E_{\lambda,2})=0$.

Finally, specialize to $\Omega=B_{r_0}$. The incident fields are divergence-free in the ordinary sense and have zero normal component on the spherical boundary. Thus $p_{\lambda,E}=p_{\lambda,H}=0$. By uniqueness, the two electric fields just constructed coincide with \eqref{piecewiseE} and \eqref{piecewiseH}, respectively. The identity $u=\tilde j_\lambda-E_{\lambda,1}$ now gives \eqref{casesR}. To identify the corresponding magnetic field explicitly, use $H_\lambda=(\mathrm{i}\lambda)^{-1}\mu^{-1}\curl E_{\lambda,2}$. Applying $\curl N_{\ell m}(k_jr)=k_jM_{\ell m}(k_jr)$ and $k_j/(\lambda\mu_j)=\eta_j^{-1}$ in each medium gives
\begin{equation*}
H_\lambda(Y_n)=\frac{k_0}{\mathrm{i}}\sqrt{\frac{2}{\pi\eta_0\ell(\ell+1)}}
\begin{cases}
\displaystyle\frac{\eta_0}{\eta_1}t_{\ell}^{2}(\lambda)M_{\ell m}^{0}(k_1r),&r<r_0,\\[1ex]
M_{\ell m}^{0}(k_0r)+A_\ell^{2}(\lambda)M_{\ell m}^{-}(k_0r),&r>r_0.
\end{cases}
\end{equation*}
The factor $\eta_0/\eta_1$ comes from factoring out the exterior normalization $1/\eta_0$; the interior magnetic field carries $1/\eta_1$. Conversely, applying $\mathrm{i}\epsilon^{-1}\curl/\lambda$ to this formula recovers \eqref{piecewiseH}.
\end{proof}

\begin{proposition}\label{outgoingsol}
For $\Omega=B_{r_0}$, the matrices $A_{\ell}^j(\lambda)$, $j=1,2$ coincide exactly with the T-matrix entries. The scattering amplitudes are obtained from these outgoing expansion coefficients by the phase and normalization factors displayed below.
\end{proposition}
\begin{proof}
Insert $h_\ell^{(1)}(y)=(-\mathrm{i})^{\ell+1}\mathrm{e}^{\mathrm{i}y}y^{-1}(1+O(y^{-1}))$ from \eqref{eq:largexi} into the Hansen fields. The recurrence \eqref{recurrence} gives $(yh_\ell^{(1)}(y))'/y=(-\mathrm{i})^\ell\mathrm{e}^{\mathrm{i}y}y^{-1}(1+O(y^{-1}))$, while the radial component of $N_{\ell m}^-$ is $O(\mathrm{e}^{\mathrm{i}k_0r}r^{-2})$. Applying these asymptotics to \eqref{piecewiseE}--\eqref{piecewiseH} gives, for fixed $\ell$ and $k_0$ in a compact subset of $\mathbb C\setminus\{0\}$,
\begin{align}
E_{\lambda,1}(Y_n)-\tilde{j}_{\lambda}(Y_n)=(-\mathrm{i})^{\ell+1}\sqrt{\frac{2\eta_0}{\pi}}A_{\ell}^1(\lambda)\frac{\mathrm{e}^{\mathrm{i}k_0r}}{r}\Psi_{1\ell m}(1+O(r^{-1}))
\end{align}
and
\begin{align}
E_{\lambda,2}(Y_n)-\frac{1}{k_0}\curl(\tilde{j}_{\lambda}(Y_n))=(-\mathrm{i})^{\ell}\sqrt{\frac{2\eta_0}{\pi}}A_{\ell}^2(\lambda)\frac{\mathrm{e}^{\mathrm{i}k_0r}}{r}(\Psi_{2\ell m}+O(r^{-1})).
\end{align}
\end{proof}

The following assertion holds for every bounded Lipschitz inclusion.
\begin{proposition}\label{nohalfplane}
There are no nonzero resonances for $P_E$ or $P_H$ with $\Im\lambda\geq0$.
\end{proposition}
\begin{proof}
An outgoing solution with $\Im\lambda>0$ decays exponentially and, by the transmission conditions, belongs to $\mathcal D(P_E)$. It would therefore be an $L^2$ eigenfunction with eigenvalue $\lambda^2\notin[0,\infty)$, contradicting non-negativity and self-adjointness. For real $\lambda\ne0$, Theorem~\ref{Rellich} gives uniqueness. The magnetic case is identical.
\end{proof}

\begin{proposition}\label{noupper}
The poles of the scattering matrix are in the lower half-plane. For the ball, set $z=x_0=k_0r_0$, so that $x_1=m_0z$.
The TE resonance occurs when
\begin{align}\label{TEcond}
\frac{\psi_{\ell}'(m_0z)}{\psi_{\ell}(m_0z)}=
\frac{\eta_1}{\eta_0}\frac{\xi_{\ell}'(z)}{\xi_{\ell}(z)}.
\end{align}
The TM resonance occurs when
\begin{align}\label{TMcond}
\frac{\psi_{\ell}'(m_0z)}{\psi_{\ell}(m_0z)}=\frac{\eta_0}{\eta_1}\frac{\xi_{\ell}'(z)}{\xi_{\ell}(z)}.
\end{align}
\end{proposition}
\begin{proof}
A purely outgoing TE field has the form $tM_{\ell m}^0(k_1r)$ inside and $cM_{\ell m}^-(k_0r)$ outside. The transmission conditions are
\begin{equation*}
\begin{pmatrix}
j_\ell(x_1)&-h_\ell^{(1)}(x_0)\\
\mu_1^{-1}\psi_\ell'(x_1)&-\mu_0^{-1}\xi_\ell'(x_0)
\end{pmatrix}
\begin{pmatrix}t\\c\end{pmatrix}=0.
\end{equation*}
This system has a nonzero solution precisely when $D_\ell^1(\lambda)=0$, which gives \eqref{TEcond} when $\psi_\ell(m_0z)\xi_\ell(z)\ne0$. Interchanging $\mu$ and $\epsilon$ gives the TM condition. Proposition~\ref{nohalfplane} locates the poles in the lower half-plane.
\end{proof}

\section{Proof of Theorem \texorpdfstring{\ref{maintheorem}}{1.4}}\label{sec:proof1}

\begin{proof}[Proof of Theorem \ref{maintheorem}]
The spectral representation below follows from Proposition~\ref{spectralresolution}.

We set for $\Im\lambda>0$
\begin{align}
G_{\lambda,E}(x,y)=\sum_n\int\limits_0^{\infty}\frac{E_{\sigma}(Y_n)(x)\otimes \overline{E_{\sigma}(Y_n)(y)}}{\sigma^2-\lambda^2}\,\mathrm{d}\sigma.
\end{align}
It follows that the resolvent can be interpreted as
\begin{align}
(R_{\lambda,E}f)(x)=\int\limits_{\mathbb{R}^3}G_{\lambda,E}(x,y)f(y)\epsilon(y)\,\mathrm{d}y \quad \div(\epsilon f)=0.
\end{align}
Let $E_\lambda:L^2_t(\mathbb{S}^2;\mathbb{C}^3)\to L^2_{\loc}$, $Y\mapsto E_\lambda(Y)$, be as in Definition~\ref{def:Elam}, and let $(\chi E_{\bar\lambda})^*$ be the adjoint of $\chi E_{\bar\lambda}$ with respect to $L^2(\epsilon\,\mathrm{d}x)$ and $L^2(\mathbb{S}^2)$, so that $(\chi E_\lambda)(\chi E_{\bar\lambda})^*f=\sum_n\chi E_\lambda(Y_n)\langle \chi f,E_{\bar\lambda}(Y_n)\rangle_{(L^2(\epsilon\,\mathrm{d}x))^3}$. Since $\lambda\mapsto\overline{E_{\bar\lambda}(Y)(x)}$ is holomorphic wherever $E_{\bar\lambda}$ is, so is $(\chi E_{\bar\lambda})^*$. After spatial localization, its mode expansion converges locally uniformly in operator norm away from zero and the zeros of $D_\ell^j$, by the fixed-argument large-order Bessel estimates: for fixed $\lambda$ and large $\ell$, we have $|D_{\ell}^j|\geq c|\psi_{\ell}(x_1)\xi_{\ell}(x_0)|$ which controls the matrix coefficients while $j_{\ell}(kr)=O(\frac{(|kR|)^{\ell}}{(2\ell+1)!!})$. For real $\chi$, the jump identity in Proposition~\ref{spectralresolution} gives

\begin{equation*}
\chi\tilde R_{\lambda,E}\chi
=\chi\tilde R_{-\lambda,E}\chi
+\frac{\mathrm{i}\pi}{\lambda}(\chi E_\lambda)(\chi E_{\bar\lambda})^*.
\end{equation*}

For $\Im\lambda<0$, both $\tilde R_{-\lambda,E}$ and $(\chi E_{\bar\lambda})^*$ are holomorphic. Thus the continuation has possible poles only at zeros of $D_\ell^j$. The gradient term $-\lambda^{-2}(I-\Pi_E)$ has no jump; the same argument applies to $P_H$.
Here the jump identity is used in the form for $\tilde R_{\pm\lambda,E}$ and $f\in L^2_{\compp}$ stated in Proposition~\ref{spectralresolution}. Both sides are holomorphic in their half planes and their boundary values agree on $(C_1,\infty)$, the boundary values of $\chi\tilde R_{\lambda,E}\chi$ from $\Im\lambda>0$ being continuous (mode by mode, together with the uniform bound on the tail of the mode sum). Hence they define one holomorphic function across the real axis.
 By Propositions~\ref{nohalfplane} and~\ref{noupper}, it remains to exclude zeros of the transmission determinants in the lower half-plane.  The denominators were originally defined as functions of the physical frequency $\lambda$. Therefore we set
\begin{align}
\mathcal{D}_{\ell}^j(z)=D_{\ell}^j\left(\frac{z}{r_0\sqrt{\epsilon_0\mu_0}}\right) \quad z=x_0=k_0r_0\quad x_1=m_0z.
\end{align}
Thus
\begin{align}
&\mathcal{D}_{\ell}^1(z)=\eta_0\xi_{\ell}(z)\psi_{\ell}'(m_0z)-\eta_1\psi_{\ell}(m_0z)\xi_{\ell}'(z)\\&
\mathcal{D}_{\ell}^2(z)=\eta_0\psi_{\ell}(m_0z)\xi_{\ell}'
(z)-\eta_1\xi_{\ell}(z)\psi_{\ell}'(m_0z).
\end{align}
First we note that
\begin{align}
\frac{\psi_{\ell}'(z)}{\psi_{\ell}(z)}=\frac{1}{2z}+\frac{J'_{\nu}(z)}{J_{\nu}(z)} \quad \textrm{and} \quad
\frac{\xi_{\ell}'(z)}{\xi_{\ell}(z)}=\frac{1}{2z}+\frac{(H^{(1)}_{\nu})'(z)}{H^{(1)}_{\nu}(z)}.
\end{align}
Assume first $0<m_0<1$; the endpoint $m_0=1$ is treated at the end of the proof. Let $M,K$ be large constants. Throughout the proof, $\nu_0$ is a natural number depending on $M,K$, and we assume that $r_0\sqrt{\epsilon_0\mu_0}\,C_1\geq K\nu_0^2/m_0$ (this is the condition $\Re z\geq K\nu_0^2/m_0$ in the variable $z=r_0\sqrt{\epsilon_0\mu_0}\,\lambda$). This implies that for all $\nu, \nu<\nu_0$ only region 4 happens. In the rest of the proof $C_2$ denotes the width of the strip in the variable $z$, chosen smaller than both \eqref{eq:C2threshold} and the Region~4 bound; the constant of Theorem~\ref{maintheorem} is then $C_2/(r_0\sqrt{\epsilon_0\mu_0})$. We use the symmetry $\mathcal{D}_{\ell}^j(-\overline{z})=-\overline{\mathcal{D}_{\ell}^j(z)}$ so that we only consider $\mathrm{Re}z>0$. It therefore suffices to consider $-C_2\leq\Im z<0$.  We consider the regions:
\begin{enumerate}
\item[] Region 1: $|\Re z|\leq\nu-M\nu^{\frac{1}{3}}$

 \item[] Region 2: $|\Re z-\nu|<M\nu^{\frac{1}{3}}$

\item[] Region 3: $\Re z-\nu \geq M\nu^{\frac{1}{3}}$, $|\Re z|<K\frac{\nu^2}{m_0}$

\item[] Region 4: $|\Re z|\geq K\frac{\nu^2}{m_0}$
\end{enumerate}

\textbf{For region 1}, both $z$ and $m_0z$ lie in the subturning range. Set
\begin{align}
\beta_1(z)=\frac{\sqrt{\nu^2-m_0^2z^2}}{m_0z} \qquad \beta_0(z)=\frac{\sqrt{\nu^2-z^2}}{z}.
\end{align}
The uniform estimates in Appendix~\ref{app:subturning}, together with \eqref{eq:riccati-definitions}, give
\begin{align}
\frac{\psi_{\ell}'(m_0z)}{\psi_{\ell}(m_0z)}=\frac{1}{2m_0z}+\beta_1(z)(1+O(M^{-\frac{3}{2}}))
\end{align}
and
\begin{align}
\frac{\xi_{\ell}'(z)}{\xi_{\ell}(z)}=\frac{1}{2z}-\beta_0(z)(1+O(M^{-\frac{3}{2}})).
\end{align}
The TE condition \eqref{TEcond} would therefore force
\begin{align}
\beta_1(1+O(M^{-\frac{3}{2}}))+\frac{\eta_1}{\eta_0}\beta_0(1+O(M^{-\frac{3}{2}}))=\left(\frac{\eta_1}{\eta_0}-\frac{1}{m_0}\right)\frac{1}{2z}.
\end{align}
The fact that $|\Im z|\leq C_2$ and $\Re z\geq C_1$ implies that $|\arg(\beta_j)|\leq \frac{\pi}{4}$ for all $\nu\geq \nu_0$, then the real part of the left hand side is bigger than $\frac{1}{2}\Re \beta_1$.  Because $\mathrm{Re}(\beta_1)\geq \frac{\sqrt{1-m_0^2}}{m_0}>0$ and $\mathrm{Re}\beta_0\geq CM^{\frac{1}{2}}\nu^{-\frac{1}{3}}>0$.
The right hand side is sufficiently smaller provided $\Re(z)$ is chosen larger than a constant depending on $m_0,\eta_1,\eta_0$ and $M$. The resonance condition cannot be true for any $\nu\geq \nu_0$, for some $\nu_0>0$, i.e.
 $\mathcal{D}_{\ell}^{1}(z)\neq 0$. The same computation with $\eta_1$ and $\eta_0$ swapped gives that $\mathcal{D}_{\ell}^{2}(z)\neq 0$.

\textbf{For region 2}, $|\Re z-\nu|<M\nu^{1/3}$ and, for $\nu\geq\nu_0(m_0,M)$,
$m_0\Re z\leq m_0(\nu+M\nu^{1/3})\leq(1+m_0)\nu/2$.
Thus the interior argument stays uniformly below its turning point, and the subturning estimates give
\begin{align}
\frac{\psi_{\ell}'(m_0z)}{\psi_{\ell}(m_0z)}=\frac{1}{2m_0z}+\beta_1(z)(1+O(\nu^{-1}))=\frac{\sqrt{1-m_0^2}}{m_0}+O(M\nu^{-\frac{2}{3}})
\end{align}
and
\begin{align}
\mathrm{Re}\frac{\psi_{\ell}'(m_0z)}{\psi_{\ell}(m_0z)}\geq \frac{1}{2}\frac{\sqrt{1-m_0^2}}{m_0}>0 \quad \nu\geq \nu_0.
\end{align}
This is where the hypothesis that $m_0<1$ is used.
With $w=z$ and the notation of Proposition~\ref{prop:two-term-airy}, so that $\tau=(z-\nu)\nu^{-1/3}$, the additive expansions give
\begin{equation*}
\frac{\xi_\ell'(z)}{\xi_\ell(z)}
 =\nu^{-1/3}\frac{p\widehat A'+O_M(\nu^{-2/3})}
                            {\widehat A+O_M(\nu^{-2/3})}.
\end{equation*}
The compact turning set may be enlarged from $M$ to $M+1$ to include
the small imaginary part of $\tau$. For real $\tau$, $\sigma$ lies on the rays of arguments $-\pi/3$ and $2\pi/3$, which avoid the zeros of $\operatorname{Ai}$. By continuity, $|\widehat A|\geq c(M)>0$ on a sufficiently small complex neighbourhood of the compact real turning interval. Since $|\Im\tau|\leq C_2\nu^{-1/3}$, it follows that
\begin{align}
\left|\frac{\xi_{\ell}'(z)}{\xi_{\ell}(z)}\right|\leq C\nu^{-\frac{1}{3}}.
\end{align}
These can also not be equal to each other for sufficiently large $\Re z$. Sufficiently large $\Re z$ implies sufficiently large $\nu$.

\textbf{For region 3}, There are all three sub-possibilities as, $m_0\Re z<\nu-M\nu^{\frac{1}{3}}$, $|m_0\Re z-\nu|<M\nu^{\frac{1}{3}}$, and $m_0\Re z-\nu>M\nu^{\frac{1}{3}}$ are all possible.
  We set $q=\frac{z}{\nu}$, it then follows that $\mathrm{Re}q=\frac{\mathrm{Re z}}{\nu}\geq 1+M\nu^{-\frac{2}{3}}$. We also have that $|\Im z|\leq C_2$
\begin{align}
|\Im q|\leq \frac{|\Im z|}{\nu}\leq \frac{C_{2}}{\nu}\le C_{2}.
\end{align}
Therefore $q$ lies in the appropriate sector, $|\mathrm{arg}(q)|\leq \pi-\delta$ for sufficiently large $\nu$.  The asymptotics in the Appendix \ref{largeturning} hold uniformly for the analytic continuation in $q$ to this region of the complex plane. It follows that $\zeta(q)=\zeta(\mathrm{Re}q)+O(C_2/\nu)$. With $\zeta(q)$ as in \eqref{eq:zeta-oscillatory}, set
\begin{align}
W_{\nu}(\Re(q))=\mathrm{e}^{\frac{2\pi \mathrm{i}}{3}}\nu^{\frac{2}{3}}\zeta(\Re(q)), \qquad \Gamma(\Re q)=\frac{1}{\Re q}\sqrt{\frac{(\Re q)^2-1}{-\zeta(\Re q)}}.
\end{align}
Uniformly in Region~3,
\begin{align}\label{eq:exterior-logderivative}
\frac{\xi_{\ell}'(z)}{\xi_{\ell}(z)}=\nu^{-\frac{1}{3}}\mathrm{e}^{-\frac{\pi \mathrm{i}}{3}}\Gamma(q)\frac{Ai'(W_{\nu}(q))}{Ai(W_{\nu}(q))}+\frac{1}{2\nu q}+O(\nu^{-1}).
\end{align}
Since $\zeta(\Re(q))<0$ on the real ray and $\zeta$ is conformal there, $|\arg W_{\nu}(\Re(q))+\pi/3|\le\delta$
for $\nu\ge\nu_{0}$; in particular $W_{\nu}(\Re(q))$ stays in a closed sector disjoint from the zeros of
$\mathrm{Ai}$ which are on the negative real axis, and $|W_{\nu}(\Re(q))|\ge c>0$. This is because $\mathrm{Re}q-1\geq M\nu^{-\frac{2}{3}}$ in this region.
Along the real part of $q$ we have $W_{\nu}(\Re(q))$
\begin{align}
\frac{Ai'(W_{\nu}(q))}{Ai(W_{\nu}(q))}&=-(W_{\nu}(q))^{\frac{1}{2}}(1+O(|W|^{-\frac{3}{2}})),\nonumber\\
W_\nu(\Re q)&=\mathrm{e}^{-\frac{\pi \mathrm{i}}{3}}\nu^{\frac{2}{3}}|\zeta(\Re q)|.
\end{align}
 This expansion gives that
\begin{align}\label{lhsairy}
\frac{\xi_{\ell}'(z)}{\xi_{\ell}(z)}=\mathrm{i}\sqrt{1-q^{-2}}(1+O(\nu^{-1}(-\zeta)^{-\frac{3}{2}}))+O(\nu^{-1}).
\end{align}
At the lower edge since $-\zeta\sim M\nu^{-\frac{2}{3}}$ there, the relative error $O(\nu^{-1}(-\zeta)^{-\frac{3}{2}})=O(M^{-\frac{3}{2}})$.
In the first subcase,
\begin{align}
\frac{\psi_{\ell}'(m_0z)}{\psi_{\ell}(m_0z)}=\frac{1}{2m_0z}+\beta_1(z)(1+O(M^{-\frac{3}{2}})).
\end{align}
One then has that $\nu^2-m_0^2(\Re z)^2\geq cM\nu^{\frac{4}{3}}$. Therefore $\Re \beta_1\geq c'M^{\frac{1}{2}}\nu^{-\frac{1}{3}}$. On the other side, by \eqref{lhsairy}, $\Re\big(\xi_\ell'(z)/\xi_\ell(z)\big)=-\Im\gamma_0+O\big(|\gamma_0|\nu^{-1}(-\zeta)^{-\frac32}\big)+O(\nu^{-1})$ with $\gamma_0=\sqrt{1-q^{-2}}$. Since $|\gamma_0|\simeq(q-1)^{\frac12}$ and $-\zeta\simeq q-1$ for $q$ near $1$, the middle term is $\lesssim\nu^{-1}(q-1)^{-1}\leq M^{-1}\nu^{-\frac13}$ on region 3, and $|\Im\gamma_0|\lesssim C_2\nu^{-1}|\gamma_0|^{-1}\leq C_2M^{-\frac12}\nu^{-\frac23}$. Hence $|\Re(\xi_\ell'/\xi_\ell)(z)|\leq C(M^{-1}\nu^{-\frac13}+\nu^{-\frac23})$, which for $M$ and $\nu$ large is smaller than $\frac{\eta_0}{\eta_1}\cdot\frac{c'}{2}M^{\frac12}\nu^{-\frac13}$. So the real part of the left hand side of \eqref{TEcond} exceeds that of the right hand side (and likewise for \eqref{TMcond}).

For the second subcase, let $|m_0\Re z-\nu|<M\nu^{\frac13}$. Write $m_0z=\nu+\tau\nu^{\frac13}$, so $|\Re\tau|<M$ and $|\Im\tau|\le m_0C_2\nu^{-\frac13}$. Use $s,A,A'$ from Proposition~\ref{prop:two-term-airy} with $w=m_0z$. The exterior variable satisfies $\Re q\ge m_0^{-1}(1-M\nu^{-\frac23})\ge\frac12(1+m_0^{-1})>1$, so it is uniformly bounded below in this oscillatory regime. Since $\xi_\ell(z)\ne0$,

\begin{equation*}
L_+(z):=\frac{\xi'_\ell(z)}{\xi_\ell(z)}=\mathrm{i}\gamma_0(z)+O(\nu^{-1}),\qquad \gamma_0(z)=\sqrt{1-m_0^2}+O(M\nu^{-\frac23}).
\end{equation*}

Divide $\mathcal D^1_\ell(z)$ by $-\sqrt\pi\,2^{\frac16}\nu^{-\frac16}\xi_\ell(z)$ and insert the expansions from the Appendix \eqref{eq:airy-psi}, \eqref{eq:airy-psiprime}, \eqref{eq:airy-riccati-errors}. These are additive, so no division by $\psi_\ell(m_0z)$ occurs. The result is
\begin{align}\label{eq:sub2}
\frac{-\mathcal{D}^1_\ell(z)}{\sqrt\pi\,2^{\frac16}\nu^{-\frac16}\xi_\ell(z)}&=\eta_0A'+\mathrm{i}\eta_1\gamma_0\,2^{-\frac13}\nu^{\frac13}A+\mathcal F,\nonumber\\
|\mathcal F|&\le C_M\nu^{-\frac13}\big(|A|+|A'|\big),
\end{align}
where we used that $|A|+|A'|\ge c_M>0$ on the compact set $|s|\le2^{\frac13}(M+1)$ ($\mathrm{Ai}$ and $\mathrm{Ai}'$ have no common zero). We argue by contradiction that $\mathcal{D}^1_\ell(z)=0$. Then $\nu^{\frac13}|A|\le C(|A'|+|\mathcal F|)$, so $|A|\le C_M\nu^{-\frac13}$ and $|A'|\ge c_M/2$. Thus $s$ lies within $C_M\nu^{-\frac13}$ of a zero $a_k<0$ of $\mathrm{Ai}$. Since $\mathrm{Ai}''(a_k)=a_k\mathrm{Ai}(a_k)=0$, we have $A/A'=(s-a_k)(1+O(|s-a_k|^2))$. Solving \eqref{eq:sub2} for $A/A'$ gives

\begin{equation*}
s-a_k=\frac{\mathrm{i}\,2^{\frac13}\eta_0}{\eta_1\gamma_0\,\nu^{\frac13}}\big(1+O_M(\nu^{-\frac13})\big),
\end{equation*}
so that
\begin{equation*}
\Im s=\frac{2^{\frac13}\eta_0}{\eta_1\sqrt{1-m_0^2}}\,\nu^{-\frac13}\big(1+O_M(\nu^{-\frac13})\big).
\end{equation*}

On the other hand $\Im s=-2^{\frac13}m_0\,\Im z\,\nu^{-\frac13}\le2^{\frac13}m_0C_2\nu^{-\frac13}$. This is a contradiction for $\nu\ge\nu_0(M)$ as soon as
\begin{align}\label{eq:C2sub2}
C_2<\frac{\eta_0}{2m_0\eta_1\sqrt{1-m_0^2}}\quad(\text{TE}),\qquad C_2<\frac{\eta_1}{2m_0\eta_0\sqrt{1-m_0^2}}\quad(\text{TM}).
\end{align}
Both conditions follow from \eqref{eq:C2threshold}, because $\frac{\eta_0\eta_1}{(\eta_0+\eta_1)^2}\le\min(\frac{\eta_0}{\eta_1},\frac{\eta_1}{\eta_0})$ and $\sqrt{1-m_0^2}\le(1-m_0^2)^{-1/2}$.

\textbf{For the third subcase}, $m_0\Re z\geq \nu+M\nu^{\frac{1}{3}}$. Here $\psi_{\ell}(m_0z)$ is oscillatory and has zeros, so $\psi_{\ell}'(m_0z)/\psi_{\ell}(m_0z)$ is an unbounded function of $z$ and cannot be compared termwise with \eqref{lhsairy}. Throughout this subcase set, with principal branches, we set
\begin{align}\label{eq:gammadef}
\gamma_0(z)=\sqrt{1-\frac{\nu^2}{z^2}}=\sqrt{1-q^{-2}}, \qquad
\gamma_1(z)=\sqrt{1-\frac{\nu^2}{m_0^2z^2}},
\end{align}
so that $\mathrm{i}\gamma_0$ is precisely the leading term of \eqref{lhsairy}. Because $m_0<1$, the subcase hypothesis $m_0\Re z>\nu$ \emph{already} forces $\Re z>\nu/m_0$: the exterior variable is automatically far beyond its own turning point, and we do not sit near the lower edge of region 3. Consequently, for $\nu\geq \nu_0$ and $|\Im z|\leq C_2$,
\begin{align}\label{eq:gamma0bound}
\Re \gamma_0(z)\geq \sqrt{1-m_0^2}+O(\nu^{-1}), \qquad |\gamma_0(z)|\leq 1+O(\nu^{-1}),
\end{align}
while $(m_0\Re z)^2-\nu^2\geq 2M\nu^{\frac{4}{3}}$ gives
\begin{align}\label{eq:gamma1bound}
c\,M^{\frac{1}{2}}\nu^{-\frac{1}{3}}\leq \Re\gamma_1(z)\leq|\gamma_1(z)|\leq 1+O(\nu^{-1}).
\end{align}
Put $x=\Re z$, $y=\Im z$ and $g_j=\gamma_j(x)>0$. Since $|\gamma_j'|\leq C/(\nu g_j)$ along $x+\mathrm{i}s$, $|s|\leq C_2$, and $g_1^2\geq cM\nu^{-2/3}$, $g_0\geq\sqrt{1-m_0^2}$, we have uniformly
\begin{equation*}
\gamma_j(x+\mathrm{i}s)=g_j(1+O(\delta_\nu)),\qquad
\delta_\nu=\frac{1+C_2}{M\nu^{1/3}},\qquad j=0,1.
\end{equation*}

For fixed $M$, increase $\nu_0$ so that $\delta_\nu\leq M^{-3/2}$. We use $g_j$ in the real inequalities below. We write $\varepsilon(w)=\nu^2((\Re w)^2-\nu^2)^{-\frac{3}{2}}$ for the relative error in the Debye expansions in the oscillatory regime, cf.\ \S 10.19(ii) of \cite{NIST} and pp.~419--425 of \cite{olverbook}, which by the uniform complex-argument expansions in \S10.20 of \cite{NIST} are valid uniformly for $|\arg w|\leq \pi-\delta$. In this subcase $\varepsilon(m_0z)=O(M^{-\frac{3}{2}})$ and $\varepsilon(z)=O(\nu^{-1})$.

\textit{Step 1: an exact reflection identity.} Put $\xi_{\ell}^{\pm}(w)=wh_{\ell}^{(1,2)}(w)$, so that $\xi_{\ell}^{+}=\xi_{\ell}$ and $\psi_{\ell}=\frac{1}{2}(\xi_{\ell}^{+}+\xi_{\ell}^{-})$, and set
\begin{align}
L_{\pm}(w)&=\frac{(\xi^{\pm}_{\ell})'(w)}{\xi^{\pm}_{\ell}(w)}=\frac{1}{2w}+\frac{(H^{(1,2)}_{\nu})'(w)}{H^{(1,2)}_{\nu}(w)},\nonumber\\
\mathcal{R}(w)&=\frac{\xi^{-}_{\ell}(w)}{\xi^{+}_{\ell}(w)}=\frac{H^{(2)}_{\nu}(w)}{H^{(1)}_{\nu}(w)}.
\end{align}
Dividing $\mathcal D^1_\ell(z)$ by $\frac{1}{2}\xi_{\ell}(z)\xi^{+}_{\ell}(m_0z)$, which is non-zero because the Debye expansions have relative error $<\frac{1}{2}$ there, gives the identity
\begin{align}\label{eq:exactid}
\frac{2\mathcal{D}_{\ell}^1(z)}{\xi_{\ell}(z)\xi^{+}_{\ell}(m_0z)}
&=\Big[\eta_0L_{+}(m_0z)-\eta_1L_{+}(z)\Big]+\mathcal{R}(m_0z)\Big[\eta_0L_{-}(m_0z)-\eta_1L_{+}(z)\Big].
\end{align}
Hence, as long as the second bracket does not vanish,
\begin{align}\label{eq:fresnel}
\mathcal{D}_{\ell}^1(z)=0 \iff \mathcal{R}(m_0z)=\mathcal{F}(z):=\frac{\eta_1L_{+}(z)-\eta_0L_{+}(m_0z)}{\eta_0L_{-}(m_0z)-\eta_1L_{+}(z)}.
\end{align}
The point of \eqref{eq:fresnel} is that the two sides are compared in \emph{modulus}: the left side is forced to be almost unimodular by the strip, and the right side is a Fresnel reflection coefficient, which is strictly smaller than $1$ exactly because $m_0<1$.

\textit{Step 2: the interior Hankel quotient is almost unimodular.} For real $w>0$ and real order we have $H^{(2)}_{\nu}(w)=\overline{H^{(1)}_{\nu}(w)}$, so $|\mathcal{R}(w)|=1$ \emph{exactly} on the real axis. Since $\partial_w\log\mathcal{R}=L_{-}-L_{+}=-2\mathrm{i}\gamma(w)(1+O(\varepsilon(w)))+O(|w|^{-1})$ with $\gamma(w)=\sqrt{1-\nu^2w^{-2}}$, integrating along the vertical segment from $m_0\Re z$ to $m_0z$ gives
\begin{align}
\begin{split}
\log|\mathcal{R}(m_0z)|&=\Re\int\limits_0^{\Im z}\mathrm{i}m_0(L_{-}-L_{+})\big(m_0(\Re z+\mathrm{i}s)\big)\,\mathrm{d}s\\
&=2m_0\int_0^y\Re\gamma_1(x+\mathrm{i}s)\,\mathrm{d}s+O(C_2g_1M^{-3/2}+\nu^{-1}),
\end{split}
\end{align}
so that
\begin{align}\label{eq:Rlower}
\begin{gathered}
\big|\log|\mathcal{R}(m_0z)|\big|\leq 2m_0g_1|y|(1+CM^{-3/2})+C\nu^{-1},\\
|\mathcal{R}(m_0z)|\geq1-2m_0C_2g_1(1+CM^{-3/2})-C\nu^{-1}.
\end{gathered}
\end{align}
The factor $g_1$ in \eqref{eq:Rlower} is the whole point and must not be discarded: it is what compensates the degeneration of the right-hand side of \eqref{eq:fresnel} as the interior ray becomes grazing.

\textit{Step 3: the right-hand side of \eqref{eq:fresnel} is strictly sub-unimodular.} The Debye expansions give
\begin{align}
L_{+}(z)=\mathrm{i}\gamma_0+O(\nu^{-1}), \qquad
L_{\pm}(m_0z)=\pm \mathrm{i}\gamma_1+O\big(g_1M^{-\frac{3}{2}}\big)+O(\nu^{-1}).
\end{align}
Put $a=\eta_0\gamma_1$ and $b=\eta_1\gamma_0$. The denominator in \eqref{eq:fresnel} equals $-\mathrm{i}(a+b)+O(g_1M^{-\frac{3}{2}}+\nu^{-1})$ and, by \eqref{eq:gamma0bound}, has modulus at least $\frac{1}{2}\eta_1\sqrt{1-m_0^2}>0$; in particular the second bracket of \eqref{eq:exactid} does not vanish, and
\begin{align}\label{eq:Fupper}
|\mathcal{F}(z)|\leq \left|\frac{b-a}{b+a}\right|+C\big(g_1M^{-\frac{3}{2}}+\nu^{-1}\big).
\end{align}
From $|b+a|^2-|b-a|^2=4\Re(a\overline{b})$, together with $|a|+|b|\leq(\eta_0+\eta_1)(1+C\delta_\nu)$ and $\Re(a\overline b)\geq\eta_0\eta_1g_1g_0(1-C\delta_\nu)$, we obtain, after increasing $\nu_0$,
\begin{align}\label{eq:fresnelgap}
1-\left|\frac{b-a}{b+a}\right|
&=\frac{4\Re(a\overline{b})}{|b+a|\big(|b+a|+|b-a|\big)}
\geq \frac{2\Re(a\overline{b})}{(|a|+|b|)^2}
\geq c_{*}g_1,\nonumber\\
c_{*}&:=\frac{\eta_0\eta_1\sqrt{1-m_0^2}}{(\eta_0+\eta_1)^2}.
\end{align}
It is exactly here that the hypothesis $m_0<1$ is used: it keeps $\Re\gamma_0$ bounded away from zero even when $\gamma_1\rightarrow 0$, i.e. it rules out total internal reflection at grazing interior incidence.

\textit{Step 4: conclusion.} If $\mathcal{D}^1_{\ell}(z)=0$ in this subcase, then \eqref{eq:fresnel}, \eqref{eq:Rlower}, \eqref{eq:Fupper} and \eqref{eq:fresnelgap} give
\begin{align}
1-2m_0C_2g_1(1+CM^{-3/2})-C\nu^{-1}&\leq |\mathcal{R}(m_0z)|=|\mathcal{F}(z)|\nonumber\\
&\leq1-c_*g_1+C(g_1M^{-3/2}+\nu^{-1}).
\end{align}
Dividing by $g_1$ and using the lower bound \eqref{eq:gamma1bound} on the real axis in the terms $C\nu^{-1}/g_1$, we arrive at
\begin{align}
c_{*}\leq 2m_0C_2(1+CM^{-3/2})+CM^{-\frac{3}{2}}+CM^{-\frac{1}{2}}\nu^{-\frac{2}{3}}.
\end{align}
This is false for all $M\geq M_0(m_0,\eta_0,\eta_1)$ and $\nu\geq\nu_0(m_0,\eta_0,\eta_1,M)$ as soon as
\begin{align}\label{eq:C2threshold}
C_2<\frac{c_{*}}{4m_0}=\frac{\sqrt{1-m_0^2}}{4m_0}\frac{\eta_0\eta_1}{(\eta_0+\eta_1)^2}.
\end{align}
Hence $\mathcal{D}_{\ell}^1(z)\neq 0$ in the third subcase. Since $c_{*}$ is symmetric in $\eta_0,\eta_1$, interchanging $\eta_0$ and $\eta_1$ gives $\mathcal{D}_{\ell}^2(z)\neq 0$ under the same threshold \eqref{eq:C2threshold}.

\textbf{For region 4}, both arguments satisfy $|w|\geq K\nu^2$. The leading-order estimates \eqref{eq:largepsi}--\eqref{eq:largexip} therefore have uniform errors $O(K^{-1})$ at $w=z,m_0z$. With $\omega$ as in \eqref{eq:large-phase} and $\theta=\omega(m_0z)$, they give
\begin{align}
\mathrm{e}^{-\mathrm{i}\omega(z)}\mathcal{D}_\ell^1(z)
&=-\eta_0\sin\theta-\mathrm{i}\eta_1\cos\theta+O(K^{-1}),\label{eq:D1large}\\
\mathrm{e}^{-\mathrm{i}\omega(z)}\mathcal{D}_\ell^2(z)
&=\mathrm{i}\eta_0\cos\theta+\eta_1\sin\theta+O(K^{-1}).\label{eq:D2large}
\end{align}
For real $\theta$ each leading line has modulus at least
$\min(\eta_0,\eta_1)$. The leading line of \eqref{eq:D1large} vanishes iff $\mathrm{e}^{2\mathrm{i}\theta}=\frac{\eta_0+\eta_1}{\eta_0-\eta_1}$, and that of \eqref{eq:D2large} iff $\mathrm{e}^{2\mathrm{i}\theta}=\frac{\eta_1+\eta_0}{\eta_1-\eta_0}$; in both cases this forces
\begin{align}
\Im \theta=-y_*,\qquad y_*:=\frac{1}{2}\log \left|\frac{\eta_0+\eta_1}{\eta_0-\eta_1}\right|\in(0,\infty]
\end{align}
($y_*=\infty$ if $\eta_0=\eta_1$). Fix $y_1:=\frac12 y_*$ (any $y_1<y_*$ works; $y_1:=1$ if $y_*=\infty$). The moduli of the leading lines are $\pi$-periodic in $\Re\theta$, and the leading lines have no zeros on $\{|\Im\theta|\leq y_1\}$, so their modulus is $\geq c(y_1,\eta_0,\eta_1)>0$ there. Choosing $K$ large makes the uniform $O(K^{-1})$ remainders in \eqref{eq:D1large}--\eqref{eq:D2large} smaller than $c/2$. Since $\Im\theta=m_0\Im z$, Region 4 contributes no zeros for
\begin{align}
|\Im z|\leq \frac{y_1}{m_0}=\frac{1}{4m_0}\log \left|\frac{\eta_0+\eta_1}{\eta_0-\eta_1}\right|\quad (\eta_0\ne\eta_1),\qquad y_1=1\quad(\eta_0=\eta_1).
\end{align}
This completes the proof for $m_0<1$. \\
\textbf{The endpoint $m_0=1$}. For $\kappa\in\{\eta_1/\eta_0,\eta_0/\eta_1\}$ put
\begin{equation*}
d_{\ell,\kappa}(z)=\xi_\ell(z)\psi_\ell'(z)-\kappa\psi_\ell(z)\xi_\ell'(z),
\end{equation*}
\begin{equation*}
\mathcal D_\ell^1=\eta_0d_{\ell,\eta_1/\eta_0},\quad
\mathcal D_\ell^2=-\eta_1d_{\ell,\eta_0/\eta_1}.
\end{equation*}
If $\eta_0=\eta_1$, the Wronskian gives $d_{\ell,1}=-\mathrm{i}$, so there are no zeros. Otherwise set
\begin{equation*}
\rho=\frac{|\eta_1-\eta_0|}{\eta_1+\eta_0}\in(0,1),\qquad
0<c<\frac14\log\frac1\rho,
\end{equation*}
and work in $-c\leq\Im z\leq0$, $\Re z>0$. Fix $M$ large and then $\nu_0$ large; write $e_\nu=M^{-3/2}+M^{-1/2}\nu^{-2/3}$.

In Region~1 the subturning estimates give nonzero $\psi_\ell,\xi_\ell$ and

\begin{equation*}
\frac{d_{\ell,\kappa}}{\psi_\ell\xi_\ell}
=(1+\kappa)\frac{\sqrt{\nu^2-z^2}}{z}\bigl(1+O_\kappa(e_\nu)\bigr)\ne0.
\end{equation*}
Here the $1/(2z)$ terms are absorbed using
$|\nu^2-z^2|^{1/2}\gtrsim M^{1/2}\nu^{2/3}$.
In Region~2 put $s=-2^{1/3}(z-\nu)\nu^{-1/3}$, $A=\mathrm{Ai}(s)$ and $B=\mathrm{Bi}(s)$. The uniform Airy expansions in the appendix, with $H_\nu^{(1)}=J_\nu+\mathrm{i}Y_\nu$, yield
\begin{equation*}
d_{\ell,\kappa}
=\pi\bigl[(\kappa-1)AA'+\mathrm{i}(BA'-\kappa AB')\bigr]
+O_{M,\kappa}(\nu^{-2/3}).
\end{equation*}
For real $s$, a zero of the bracket would force $AA'=0$. The Wronskian $AB'-A'B=1/\pi$, (9.2.7) in \cite{NIST}, makes its imaginary part $-1/\pi$ when $A=0$ and $-\kappa/\pi$ when $A'=0$. Thus the bracket is bounded away from zero on the compact real turning interval. Since $|\Im s|\leq2^{1/3}c\nu^{-1/3}$, continuity and the remainder bound exclude zeros throughout this region for large $\nu$.

In Region~3 use $L_\pm,\mathcal R$ from \eqref{eq:fresnel} and put $\gamma=\sqrt{1-\nu^2/z^2}$. The Debye estimates give $L_\pm=\pm \mathrm{i}\gamma+O(|\gamma|e_\nu)$. Thus $L_--\kappa L_+=-\mathrm{i}(1+\kappa)\gamma(1+O_\kappa(e_\nu))\ne0$ for the chosen $M,\nu_0$, so the following division is valid. A zero of $d_{\ell,\kappa}$ would therefore imply
\begin{equation*}
|\mathcal R(z)|
=\left|\frac{(\kappa-1)L_+(z)}{L_-(z)-\kappa L_+(z)}\right|
\leq\rho+O_\kappa(e_\nu).
\end{equation*}
But $|\mathcal R(x)|=1$ for real $x$ and
$|\partial_y\log|\mathcal R(x+\mathrm{i}y)||\leq3$ by the same estimates. Hence $|\mathcal R(z)|\geq \mathrm{e}^{-3c}>\rho$, a contradiction after increasing $M,\nu_0$.
Region~4 applies unchanged at $m_0=1$ for this choice of $c$ and large $K$. Taking $C_1$ large enough to place all $\nu<\nu_0$ in Region~4 completes the proof, with $C_2=c/(r_0\sqrt{\epsilon_0\mu_0})$.
\end{proof}
\section{Generalized eigenfunctions and spectral resolution}\label{sec:geneig}
Throughout this section $\Omega$ may be any bounded Lipschitz open set, with coefficients \eqref{eq:lipschitz-coefficients}. We use the two families from Lemma~\ref{resolventdist}; spherical harmonics label the incident fields at infinity and do not impose a spherical interface.
\begin{definition}\label{def:Elam}
We set $a_{j,\ell,m}=\langle Y,\Psi_{j,\ell,m}\rangle$
\begin{align}
E_{\lambda}(Y)=\sum\limits_{\ell,m}(-\mathrm{i})^{\ell+1}a_{1,\ell,m}E_{\lambda,1}(Y_{1,\ell,m})+a_{2,\ell,m}(-\mathrm{i})^{\ell}E_{\lambda,2}(Y_{2,\ell,m})
\end{align}
This is the electric field generalized eigenfunction. For $Y\in L^2_t(\mathbb S^2)$ the sum converges in $H_{\loc}(\curl)$, by the local free-wave estimates, \eqref{eq:incident-gradient}, and the resolvent bounds in Lemma~\ref{resolventdist}.
\end{definition}
For smooth tangential $Y$, the free-wave asymptotics and the outgoing correction give
\begin{align}\label{incomingE}
E_\lambda(Y)=\sqrt{\frac{\eta_0}{2\pi}}\left(\frac{\mathrm{e}^{-\mathrm{i}k_0r}}rY+\frac{\mathrm{e}^{\mathrm{i}k_0r}}r\widetilde Y\right)+O(r^{-2}).
\end{align}
Here $\widetilde Y$ is smooth and tangential. 
The generalized eigenfunctions $E_{\lambda}(Y)$ can be viewed as tempered distributions in $\mathcal{S}^\prime(\mathbb{R}^3;\mathbb{C}^3)$. In particular $g\in C_0^{\infty}(\mathbb{R}_+)$ then $\int\limits_{\mathbb{R}}g(\lambda)E_{\lambda}(Y)\,\mathrm{d}\lambda$ is square integrable. Therefore the inner product
\begin{align}
\langle E_{\lambda}(Y),E_{\lambda'}(Y')\rangle_{L^2(\epsilon \mathrm{d}x)^3}
\end{align}
can be viewed as a bidistribution in $\mathcal{D}'(\mathbb{R}_+\times\mathbb{R}_+)$. The last inner product can be computed by taking the limit
\begin{align*}
&\lim\limits_{R\rightarrow\infty}\frac{1}{\lambda^2-\lambda'^2}\Big(\langle P_{E} E_{\lambda}(Y),\chi_RE_{\lambda'}(Y')\rangle_{L^2(\epsilon \mathrm{d}x)^3}-\langle \chi_R E_{\lambda}(Y), P_E E_{\lambda'}(Y')\rangle_{L^2(\epsilon \mathrm{d}x)^3}\Big)
\end{align*}
and then using Green's identity. Here $\chi_R$ denotes the indicator function of a large ball $B_R$ where $\overline\Omega\subset B_R$. 
Write $a^2=\eta_0/(2\pi)$, $k=k_0(\lambda)$, $k'=k_0(\lambda')$ and
$C=\langle Y,Y'\rangle$, $D(\lambda,\lambda')=\langle\widetilde Y_\lambda,\widetilde Y'_{\lambda'}\rangle$, with angular inner products on $\mathbb S^2$.
Green's identity at equal frequencies gives, $D(\lambda,\lambda)=C$.
Inserting \eqref{incomingE} into the boundary term gives
\begin{align*}
\langle E_\lambda(Y),\chi_RE_{\lambda'}(Y')\rangle
&=\epsilon_0a^2\frac{\sin((k-k')R)}{k-k'}(C+D(\lambda,\lambda'))
+o_{\mathcal D'}(1).
\end{align*}
Indeed, the mixed incoming--outgoing terms have phases $\pm(k+k')R$ and vanish distributionally on $\mathbb R_+^2$. The remaining cosine term has coefficient $(C-D)/(k-k')$ tends to zero distributionally. Since $\sin(Rx)/x\to\pi\delta(x)$ and
\begin{equation*}
2\pi\epsilon_0a^2\delta(k-k')
=\frac{\epsilon_0\eta_0}{\sqrt{\epsilon_0\mu_0}}\delta(\lambda-\lambda')
=\delta(\lambda-\lambda'),
\end{equation*}
we obtain the distributional identity:
\begin{align*}
&\langle E_{\lambda}(Y),E_{\lambda'}(Y')\rangle_{L^2(\epsilon \mathrm{d}x)^3}= \\&2\lambda\,\delta(\lambda^2-\lambda'^2)\langle Y,Y'\rangle_{L^2(\mathbb{S}^2)}=\delta(\lambda-\lambda')\langle Y,Y'\rangle_{L^2(\mathbb{S}^2)}, \lambda,\lambda'>0
\end{align*}

The generalized eigenfunctions can be used to describe evolution operators using the Spectral theorem.
\begin{proposition}\label{spectralresolution}
If $f\in L^2_{\compp}(\mathbb{R}^3;\mathbb{C}^3)\cap L_\epsilon(\div0)$, then for real $\lambda>0$
\begin{align}
(R_{\lambda,E}-R_{-\lambda,E})f=\frac{\mathrm{i}\pi}{\lambda}\sum\limits_{n}E_{\lambda}(Y_n)\langle f, E_{\overline{\lambda}}(Y_n)\rangle_{(L^2(\epsilon \mathrm{d}x))^3}
\end{align}
in $L^2_{loc}(\mathbb{R}^3)$. For the spectral measure $\mathrm{d}B_\lambda$ of $\sqrt{P_E}$, and compactly supported $f,g\in L_\epsilon(\div0)$, we have
\begin{align}
\langle \mathrm{d}B_{\lambda}f,g\rangle_{(L^2(\epsilon \mathrm{d}x))^3} =\chi_{[0,\infty)}(\lambda)
\sum\limits_{n}\langle f,E_{\lambda}(Y_n)\rangle_{(L^2(\epsilon \mathrm{d}x))^3} \langle E_{\overline{\lambda}}(Y_n),g\rangle_{(L^2(\epsilon \mathrm{d}x))^3} \,\mathrm{d}\lambda
\end{align}
For any bounded Borel function we then have that
\begin{align}
\langle h(P_E)f,g\rangle_{(L^2(\epsilon \mathrm{d}x))^3} &=\sum\limits_n\int\limits_0^{\infty} h(\lambda^2) \langle f,E_{\lambda}(Y_n)\rangle_{(L^2(\epsilon \mathrm{d}x))^3}\langle E_{\overline{\lambda}}(Y_n),g\rangle_{(L^2(\epsilon \mathrm{d}x))^3}\,\mathrm{d}\lambda.
\end{align}
The first identity (the resolvent jump) holds with $R_{\pm\lambda,E}$ replaced by $\tilde R_{\pm\lambda,E}$ for every $f\in L^2_{\compp}(\mathbb R^3;\mathbb C^3)$, without the divergence condition.
\end{proposition}
\begin{proof}
The incoming counterpart of \eqref{OR} gives
\begin{align}
R_{-\lambda,E}f=\sqrt{\frac{2\eta_0}{\pi}}\frac{\mathrm{e}^{-\mathrm{i}k_0r}}{r}(\Psi +O(r^{-1}))
\end{align}
for some $\Psi\in C^{\infty}(\mathbb{S}^2;\mathbb{C}^3)$. We then have that
\begin{align}
(R_{\lambda,E}-R_{-\lambda,E})f=\sqrt{\frac{2\eta_0}{\pi}}\left(-\frac{\mathrm{e}^{-\mathrm{i}k_0r}}{r}\Psi+\frac{\mathrm{e}^{\mathrm{i}k_0r}}{r}\tilde{\Psi}\right)+O(r^{-2})
\end{align}
as $r\rightarrow \infty$ for some $\tilde{\Psi}\in C^{\infty}(\mathbb{S}^2;\mathbb{C}^3)$ and
\begin{align}
(P_E-\lambda^2)(R_{\lambda,E}-R_{-\lambda,E})f=0.
\end{align}
By \eqref{incomingE}, $E_{\lambda}(-2\Psi)$ and $(R_{\lambda,E}-R_{-\lambda,E})f$ have the same incoming leading asymptotic, so their difference $w$ is purely outgoing. Since $(P_E-\lambda^2)w=0$ on all of $\mathbb{R}^3$ with the transmission conditions, Theorem~\ref{Rellich} gives $w=0$.  Integration by parts then gives that
\begin{equation}
\begin{split}
\langle f,E_{\lambda}(Y_n)\rangle_{(L^2(\epsilon \mathrm{d}x))^3}&= \langle (P_E-\lambda^2)R_{-\lambda,E}f,E_{\overline{\lambda}}(Y_n)\rangle_{(L^2(\epsilon \mathrm{d}x))^3}= \frac{2\mathrm{i}\lambda}{\pi}\langle \Psi,Y_n\rangle_{L^2(\mathbb{S}^2;\mathbb{C}^3)}
\end{split}
\end{equation}
The boundary term at $|x|=R$ is $\mu_0^{-1}\int_{\partial B_R}\big[(\hat r\times\curl u)\cdot\overline{E}-(\hat r\times\curl\overline{E})\cdot u\big]\mathrm{d}S$ with $u=R_{-\lambda,E}f$; only the incoming part of $E$ contributes in the limit, giving the result. Because the $Y_n$ are a basis for $L^2_t(\mathbb{S}^2;\mathbb{C}^3)$, and $\Psi$ is tangential, we can expand $\Psi$ as
\begin{align}
\Psi=\sum\limits_n\langle \Psi,Y_n\rangle_{L^2(\mathbb{S}^2;\mathbb{C}^3)}Y_n
\end{align}
and the sum converges in $C^{\infty}(\mathbb{S}^2;\mathbb{C}^3)$. Thus
\begin{align}
\Psi=-\frac{\mathrm{i}\pi}{2\lambda}\sum\limits_n\langle f,E_{\overline{\lambda}}(Y_n)\rangle_{(L^2(\epsilon \mathrm{d}x))^3}Y_n.
\end{align}
Consequently,
\begin{align}
(R_{\lambda,E}-R_{-\lambda,E})f=E_{\lambda}(-2\Psi)=\frac{\mathrm{i}\pi}{\lambda}\sum\limits_nE_{\lambda}(Y_n)\langle f,E_{\overline{\lambda}}(Y_n)\rangle.
\end{align}
This agrees with Stone's formula $R_{\lambda,E}-R_{-\lambda,E}=2\pi \mathrm{i}\,\delta(P_E-\lambda^2)$ for the $\delta(\lambda-\lambda')$-normalized family.

The jump formula also holds for $\tilde R_{\lambda,E}$; its gradient term has no jump. The limiting absorption argument in Lemma~\ref{resolventdist} gives continuous, locally bounded matrix elements of $R_{\lambda,E}\Pi_E$ for compactly supported $L^2$ data on compact subintervals of $(0,\infty)$. The gradient term in $\tilde R_{\lambda,E}$ is analytic there. Stone's formula therefore gives absolutely continuous spectral measures for this dense class and hence for all data in $L_\epsilon(\div0)$. On $L_\epsilon(\div0)$, $P_Eu=0$ implies $\curl u=0$, and \eqref{HD} then gives $u=0$. Thus
\begin{equation*}
\mathrm{d}B_\lambda=\frac{\lambda}{\pi \mathrm{i}}(R_{\lambda,E}-R_{-\lambda,E})\,\mathrm{d}\lambda,
\qquad \lambda>0,
\end{equation*}
accounts for the entire spectral measure of $\sqrt{P_E}$. Substitution of the jump formula proves the assertion, with extension to $L_\epsilon(\div0)$ by the spectral theorem. For $f\in L^2_{\compp}$ without the divergence condition, $u=\tilde R_{-\lambda,E}f$ solves $(\epsilon^{-1}\curl\mu^{-1}\curl-\lambda^2)u=f$, so $\div(\epsilon u)=-\lambda^{-2}\div(\epsilon f)$ is compactly supported and $u$ is an incoming solution outside $\supp f$; the argument above applies verbatim. Moreover $\langle f,E_{\overline\lambda}(Y_n)\rangle_{(L^2(\epsilon \mathrm{d}x))^3}=\langle \Pi_Ef,E_{\overline\lambda}(Y_n)\rangle_{(L^2(\epsilon \mathrm{d}x))^3}$, since $\div(\epsilon E_{\overline\lambda}(Y_n))=0$ and $f$ has compact support, and the gradient term of $\tilde R_{\pm\lambda,E}$ is even in $\lambda$.
 \end{proof}

\section{Proof of Theorem \texorpdfstring{\ref{maintheorem2}}{1.5}}
We assume that $m_0>1$ and use the notation of \eqref{eq:riccati-definitions}. We define
\begin{align}
\vartheta_{\ell}(w)=wy_{\ell}(w) \quad \xi_{\ell}=\psi_{\ell}+\mathrm{i}\vartheta_{\ell}
\end{align}
It follows that $\psi_{\ell}$ is entire and $\vartheta_{\ell}$ is holomorphic on $\mathbb{C}\setminus\{0\}$, these are both real on $(0,\infty)$ and we have that
\begin{align}\label{wgm:ode}
v''(w)+\left(1-\frac{\ell(\ell+1)}{w^2}\right)v(w)=0 \quad v=\psi_{\ell},\vartheta_{\ell}, \quad \psi_{\ell}\vartheta_{\ell}'-\psi_{\ell}'\vartheta_{\ell}=1.
\end{align}
The $\psi_{\ell}, \psi_{\ell}'$ have no common zero and neither do $\xi_{\ell}, \xi_{\ell}'$ because $\psi_{\ell}\xi_{\ell}'-\psi_{\ell}'\xi_{\ell}=\mathrm{i}$. For $\kappa>0$ we put
\begin{align}
&\mathcal{D}_{\kappa,\ell}(z)=\xi_{\ell}(z)\psi_{\ell}'(m_0z)-\kappa\psi_{\ell}(m_0z)\xi_{\ell}'(z)=\mathcal{P}_{\kappa,\ell}(z)+\mathrm{i}\mathcal{Q}_{\kappa,\ell}(z), \\&
\mathcal{P}_{\kappa,\ell}(z)=\psi_{\ell}(z)\psi_{\ell}'(m_0z)-\kappa\psi_{\ell}(m_0z)\psi_{\ell}'(z),\\&
\mathcal{Q}_{\kappa,\ell}(z)=\vartheta_{\ell}(z)\psi_{\ell}'(m_0z)-\kappa\psi_{\ell}(m_0z)\vartheta_{\ell}'(z)
\end{align}
The denominators of Section~\ref{sec:proof1} are
\begin{align}\label{wgm:DD}
\mathcal{D}^1_\ell=\eta_0\,\mathcal{D}_{\eta_1/\eta_0,\ell},\qquad \mathcal{D}^2_\ell=-\eta_1\,\mathcal{D}_{\eta_0/\eta_1,\ell},
\end{align}
and $\mathcal{P}_{\kappa,\ell},\mathcal{Q}_{\kappa,\ell}$ are real on $(0,\infty)$.
Fix $\frac{1}{m_0}<t_1<t_2<1$ and set
\begin{align}
\Lambda(t)=\mathrm{arcsech}t-\sqrt{1-t^2} \quad 0<t<1.
\end{align}
Then it follows that $\Lambda>0$ and also $\Lambda'(t)=-\frac{\sqrt{1-t^2}}{{t}}<0$ on $(0,1)$. The window $(\frac{1}{m_0},1)$ is non-empty exactly when $m_0>1$; for $z=t\nu$ with $t$ in this window the exterior argument $z$ is below its turning point (tunnelling region) while the interior argument $m_0z$ is past its turning point.

\begin{lemma}\label{wgm:debye}
Let $0<\tau_1<\tau_2<1<s_1<s_2$. We write $t=\mathrm{sech}\alpha$ and $s=\sec \beta$ and set $\Omega_{\nu}(s)=\nu(\tan\beta-\beta)-\frac{\pi}{4}$, so that $\sinh\alpha=\frac{\sqrt{1-t^2}}{t}$, $\sin\beta=\frac{\sqrt{s^2-1}}{s}$ and $\alpha-\tanh\alpha=\Lambda(t)$. As $\nu\rightarrow\infty$, uniformly for $t\in [\tau_1,\tau_2]$ and $s\in [s_1,s_2]$, we have
\begin{enumerate}
\item \begin{align}
\vartheta_{\ell}(\nu t)=-(\sinh \alpha)^{-\frac{1}{2}}\mathrm{e}^{\nu\Lambda(t)}(1+O(\nu^{-1}))
\end{align}
\item \begin{align}
\frac{\vartheta_{\ell}'(\nu t)}{\vartheta_{\ell}(\nu t)}=-\mathrm{sinh}\alpha+O(\nu^{-1})
\end{align}
\item
\begin{align}
\psi_{\ell}(\nu s)&=(\sin\beta)^{-\frac{1}{2}}(\cos \Omega_{\nu}(s)+O(\nu^{-1})),\nonumber\\
\psi_{\ell}'(\nu s)&=-(\sin \beta)^{\frac{1}{2}}(\sin \Omega_{\nu}(s)+O(\nu^{-1}))
\end{align}
\end{enumerate}
\end{lemma}
\begin{proof}
We use
\begin{align}
\vartheta_{\ell}(w)=\sqrt{\frac{\pi w}{2}}Y_{\nu}(w), \quad \vartheta_{\ell}'=\frac{\vartheta_{\ell}}{2w}+\sqrt{\frac{\pi w}{2}}Y_{\nu}',\quad \psi_{\ell}'=\frac{\psi_{\ell}}{2w}+\sqrt{\frac{\pi w}{2}}J_{\nu}'
\end{align}
The points then follow using the Debye expansions (10.19.3), (10.19.4), (10.19.6), (10.19.7) in \cite{NIST}: (i), (ii) from $Y_\nu(\nu\,\mathrm{sech}\,\alpha)$, $Y_\nu'(\nu\,\mathrm{sech}\,\alpha)$ (note $\sqrt{t/\tanh\alpha}=(\sinh\alpha)^{-1/2}$), (iii) from $J_\nu(\nu\sec\beta)$, $J_\nu'(\nu\sec\beta)$ (note $\sqrt{\sec\beta/\tan\beta}=(\sin\beta)^{-1/2}$, $\sqrt{\sec\beta\,\sin2\beta/2}=(\sin\beta)^{1/2}$); the term $\psi_\ell/(2w)$ is $O(\nu^{-1})$.
\end{proof}

\begin{lemma}\label{wgm:upper}
\begin{enumerate}
\item $|\psi_{\ell}(w)|\leq |w|\mathrm{e}^{|\Im w|}$ for all $w\in\mathbb{C}$;
\item We let $0<\tau_-<\tau_+<1$. Then there is a constant $C$ such that for all $\ell\geq 1$ and $z=x+\mathrm{i}y$ with $\frac{x}{\nu}\in [\tau_-,\tau_+]$ and $|y|\leq 3$ we have
\begin{align}
|\psi_{\ell}(z)|\leq C\nu^{\frac{1}{2}}\mathrm{e}^{-\nu\Lambda(\frac{x}{\nu})} \quad |\vartheta_{\ell}(z)|\leq C\nu^{\frac{1}{2}}\mathrm{e}^{\nu\Lambda(\frac{x}{\nu})}
\end{align}
\end{enumerate}
\end{lemma}

\begin{proof}
(i) follows from $j_\ell(w)=\frac{(-\mathrm{i})^\ell}{2}\int_{-1}^1\mathrm{e}^{\mathrm{i}ws}P_\ell(s)\,\mathrm{d}s$, see (10.54.2) in \cite{NIST}, and $|P_\ell|\le1$ on $[-1,1]$.

(ii) Let $\tau=x/\nu$ and $\alpha=\mathrm{arcsech}\,\tau$, so that $\tau\cosh\alpha=1$ and $\tau\sinh\alpha=\tanh\alpha$.

\emph{Bound for $J_\nu$.} For $\Re z>0$, $J_\nu(z)=\frac{1}{2\pi \mathrm{i}}\int_{\infty-\pi \mathrm{i}}^{\infty+\pi \mathrm{i}}\mathrm{e}^{z\sinh w-\nu w}\,\mathrm{d}w$, see \S10.9 in \cite{NIST}. The integrand is entire and decays along $\Im w=\pm\pi$, so the path may be taken to be: $\Im w=-\pi$ from $+\infty$ to $\alpha-\pi \mathrm{i}$, then the segment $\Re w=\alpha$, then $\Im w=\pi$ back to $+\infty$.
\begin{itemize}
\item On the segment $w=\alpha+\mathrm{i}\theta$: $\Re(z\sinh w-\nu w)=x\sinh\alpha\cos\theta-y\cosh\alpha\sin\theta-\nu\alpha\le-\nu\Lambda(\tau)+|y|/\tau$.
\item On the rays $w=\sigma\pm\pi \mathrm{i}$, $\sigma\ge\alpha$: $\sinh w=-\sinh\sigma$ and $\Re(z\sinh w-\nu w)=-x\sinh\sigma-\nu\sigma\le-\nu\Lambda(\tau)-\nu(\sigma-\alpha)$.
\end{itemize}
Hence $|J_\nu(z)|\le(\mathrm{e}^{3/\tau_-}+1)\mathrm{e}^{-\nu\Lambda(\tau)}$.

\emph{Bound for $Y_\nu$.} Since $\cos\nu\pi=0$, (10.9.7) in \cite{NIST} gives $Y_\nu(z)=\frac1\pi\int_0^\pi\sin(z\sin\theta-\nu\theta)\,\mathrm{d}\theta-\frac1\pi\int_0^\infty \mathrm{e}^{\nu\sigma-z\sinh\sigma}\,\mathrm{d}\sigma$. The first term is at most $\cosh y\le \mathrm{e}^3$. In the second term, $|\mathrm{e}^{\nu\sigma-z\sinh\sigma}|=\mathrm{e}^{\nu(\sigma-\tau\sinh\sigma)}\le \mathrm{e}^{\nu\Lambda(\tau)}$, the maximum being attained at $\sigma=\alpha$. Choose $T=T(\tau_-)$ with $\tau_-\sinh\sigma\ge2\sigma$ for $\sigma\ge T$; then the integrand is bounded by $\mathrm{e}^{-\nu\sigma}$ for $\sigma\ge T$. Hence $|Y_\nu(z)|\le \mathrm{e}^3+\pi^{-1}(T\mathrm{e}^{\nu\Lambda(\tau)}+\nu^{-1})\le C\mathrm{e}^{\nu\Lambda(\tau)}$. The factor $\sqrt{\pi|z|/2}\le C\nu^{\frac12}$ completes the proof.
\end{proof}

\subsection{Real quasimodes}

\begin{lemma}\label{wgm:quasi}
Let $\kappa>0$. There are $\ell_0$ and $c_1>0$ such that for each $\ell\ge\ell_0$ there is $x_\ell\in(t_1\nu,t_2\nu)$ with
\begin{align*}
\mathcal{Q}_{\kappa,\ell}(x_\ell)=0,\qquad |\mathcal{Q}_{\kappa,\ell}'(x_\ell)|\ge c_1\mathrm{e}^{\nu\Lambda(x_\ell/\nu)}.
\end{align*}
\end{lemma}

\begin{proof}
Let $I_\nu=[t_1\nu,t_2\nu]$. For $x\in I_\nu$ put $t=x/\nu$ and $s=m_0t\in[m_0t_1,m_0t_2]\subset(1,m_0)$, define $\alpha,\beta$ as in Lemma \ref{wgm:debye}, and set
\begin{align*}
&a=\sin\beta,\quad b=\sinh\alpha,\quad \mathcal{A}=\sqrt{a^2+\kappa^2b^2},\\
&\varpi\in(0,\tfrac{\pi}{2})\ \text{with}\ \mathcal{A}\cos\varpi=\kappa b,\ \mathcal{A}\sin\varpi=a,
\end{align*}
all smooth functions of $t\in[t_1,t_2]$. Let $a_1=a(t_1)>0$. By Lemma \ref{wgm:debye}(i), $\vartheta_\ell\neq0$ on $I_\nu$ for large $\nu$, so $\mathcal{Q}_{\kappa,\ell}=\vartheta_\ell\,q$ there, with
\begin{equation*}
q(x)=\psi_\ell'(m_0x)-\kappa L(x)\psi_\ell(m_0x),\qquad L=\vartheta_\ell'/\vartheta_\ell,
\end{equation*}
real valued on $I_\nu$. With $\Omega(x)=\Omega_\nu(m_0x/\nu)+\varpi(x/\nu)$, Lemma \ref{wgm:debye}(ii), (iii) give
\begin{align}\label{wgm:q}
q(x)=a^{-\frac12}\mathcal{A}\cos\Omega(x)+O(\nu^{-1}).
\end{align}
For $q'$ use \eqref{wgm:ode} and $\ell(\ell+1)=\nu^2-\frac14$:
\begin{itemize}
\item $\frac{\mathrm{d}}{\mathrm{d}x}\psi_\ell'(m_0x)=-m_0\big(a^2+O(\nu^{-2})\big)\psi_\ell(m_0x)$, because $1-\nu^2/(m_0x)^2=1-s^{-2}=a^2$;
\item $L'=\frac{\ell(\ell+1)}{x^2}-1-L^2=b^2-L^2+O(\nu^{-2})=O(\nu^{-1})$, by Lemma \ref{wgm:debye}(ii) and $t^{-2}-1=b^2$.
\end{itemize}
Hence $q'=-m_0a^2\psi_\ell(m_0x)+m_0\kappa b\,\psi_\ell'(m_0x)+O(\nu^{-1})$, and by Lemma \ref{wgm:debye}(iii)
\begin{align}\label{wgm:qprime}
q'(x)=-m_0a^{\frac12}\mathcal{A}\sin\Omega(x)+O(\nu^{-1}).
\end{align}

\emph{The phase.} Since $\frac{\mathrm{d}}{\mathrm{d}s}\big(\sqrt{s^2-1}-\mathrm{arcsec}\,s\big)=\frac{\sqrt{s^2-1}}{s}$, we have $\frac{\mathrm{d}}{\mathrm{d}x}\Omega_\nu(m_0x/\nu)=m_0a\ge m_0a_1$, while $\frac{\mathrm{d}}{\mathrm{d}x}\varpi(x/\nu)=O(\nu^{-1})$. So $\Omega$ is increasing on $I_\nu$, and $\Omega(t_2\nu)-\Omega(t_1\nu)\ge\frac12m_0a_1(t_2-t_1)\nu\ge2\pi$ for $\nu\ge\nu_0$. Choose $k\in\mathbb{Z}$ and $x_-<x_+$ in $I_\nu$ with $\Omega(x_-)=k\pi$ and $\Omega(x_+)=(k+1)\pi$.

\emph{Sign change.} As $a\le1$, $a^{-\frac12}\mathcal{A}\ge a_1$. So \eqref{wgm:q} gives $(-1)^kq(x_-)\ge a_1-C\nu^{-1}>0>(-1)^kq(x_+)$, and there is $x_\ell\in(x_-,x_+)$ with $q(x_\ell)=0$, i.e.\ $\mathcal{Q}_{\kappa,\ell}(x_\ell)=0$.

\emph{Derivative.} At $x_\ell$, \eqref{wgm:q} gives $|\cos\Omega(x_\ell)|\le C/(a_1\nu)$, hence $|\sin\Omega(x_\ell)|\ge\frac12$, and \eqref{wgm:qprime} gives $|q'(x_\ell)|\ge\frac14m_0a_1^{3/2}$. Finally $\mathcal{Q}_{\kappa,\ell}'(x_\ell)=\vartheta_\ell(x_\ell)q'(x_\ell)$, and by Lemma \ref{wgm:debye}(i), $|\vartheta_\ell(x_\ell)|\ge\frac12(\sinh\alpha(t_1))^{-\frac12}\mathrm{e}^{\nu\Lambda(x_\ell/\nu)}$. This gives the claim with $c_1=\frac18m_0a_1^{3/2}(\sinh\alpha(t_1))^{-\frac12}$.
\end{proof}

\subsection{From quasimodes to resonances}

\begin{proposition}\label{wgm:rouche}
Let $\kappa>0$ and let $x_\ell$ be as in Lemma \ref{wgm:quasi}. There are $\ell_0$ and $C>0$ such that for $\ell\ge\ell_0$, $\mathcal{D}_{\kappa,\ell}$ has exactly one zero $z_\ell$ in the disc $\{|z-x_\ell|<\rho_\ell\}$, where

\begin{equation*}
\rho_\ell=C\nu^{\frac32}\mathrm{e}^{-2\nu\Lambda(x_\ell/\nu)}.
\end{equation*}

\end{proposition}

\begin{proof}
Let $\Lambda_\ell=\Lambda(x_\ell/\nu)\ge\Lambda(t_2)>0$, and apply Lemma \ref{wgm:upper}(ii) with $\tau_-=\frac{t_1}{2}$, $\tau_+=\frac{1+t_2}{2}$. For $\nu$ large, every $z$ with $|z-x_\ell|\le2$ satisfies $\Re z/\nu\in[\tau_-,\tau_+]$, $|\Im z|\le2$, and $\nu|\Lambda(\Re z/\nu)-\Lambda_\ell|\le2\|\Lambda'\|_{L^\infty[\tau_-,\tau_+]}$.

\emph{Upper bounds.} On $\{|z-x_\ell|\le2\}$, Lemma \ref{wgm:upper} gives
$|\psi_\ell|\le C\nu^{\frac12}\mathrm{e}^{-\nu\Lambda_\ell}$, $|\vartheta_\ell|\le C\nu^{\frac12}\mathrm{e}^{\nu\Lambda_\ell}$ and $|\psi_\ell(m_0\,\cdot)|\le C\nu$. Cauchy's estimates on discs of radius $1$ give the same bounds, with new constants, for the derivatives of order $\le3$ of these three functions on $\{|z-x_\ell|\le1\}$. Since $\mathcal{P}_{\kappa,\ell}$ and $\mathcal{Q}_{\kappa,\ell}''$ are finite sums of products of such derivatives (one factor of $\psi_\ell^{(i)}$ or $\vartheta_\ell^{(i)}$, one of $(\psi_\ell(m_0\cdot))^{(j)}$),
\begin{align}\label{wgm:PQbounds}
\sup_{|z-x_\ell|\le1}|\mathcal{P}_{\kappa,\ell}(z)|\le B_2\nu^{\frac32}\mathrm{e}^{-\nu\Lambda_\ell},\qquad
\sup_{|z-x_\ell|\le1}|\mathcal{Q}_{\kappa,\ell}''(z)|\le B_3\nu^{\frac32}\mathrm{e}^{\nu\Lambda_\ell}.
\end{align}

\emph{Choice of radius.} Set $\rho_\ell=4B_2c_1^{-1}\nu^{\frac32}\mathrm{e}^{-2\nu\Lambda_\ell}$. For $\ell\ge\ell_0$ we have $\rho_\ell\le1$ and $\rho_\ell B_3\nu^{\frac32}\le c_1$, since $\Lambda_\ell\ge\Lambda(t_2)>0$.\\
\emph{Lower bound for $\mathcal{Q}_{\kappa,\ell}$.} Since $\mathcal{Q}_{\kappa,\ell}(x_\ell)=0$, Taylor's formula, \eqref{wgm:PQbounds} and Lemma \ref{wgm:quasi} give, on $|z-x_\ell|=\rho_\ell$,
\begin{align*}
|\mathcal{Q}_{\kappa,\ell}(z)-\mathcal{Q}_{\kappa,\ell}'(x_\ell)(z-x_\ell)|&\le\tfrac12\rho_\ell^2B_3\nu^{\frac32}\mathrm{e}^{\nu\Lambda_\ell}\\
&\le\tfrac12c_1\rho_\ell \mathrm{e}^{\nu\Lambda_\ell}\le\tfrac12|\mathcal{Q}_{\kappa,\ell}'(x_\ell)(z-x_\ell)|.
\end{align*}
Hence $|\mathcal{Q}_{\kappa,\ell}(z)|\ge\frac12c_1\rho_\ell \mathrm{e}^{\nu\Lambda_\ell}=2B_2\nu^{\frac32}\mathrm{e}^{-\nu\Lambda_\ell}>|\mathcal{P}_{\kappa,\ell}(z)|$ on the circle.

\emph{Rouch\'e, twice.} By the first inequality, $\mathcal{Q}_{\kappa,\ell}$ has exactly one zero in the disc, like $z\mapsto\mathcal{Q}_{\kappa,\ell}'(x_\ell)(z-x_\ell)$. By the second, $|\mathcal{D}_{\kappa,\ell}-\mathrm{i}\mathcal{Q}_{\kappa,\ell}|=|\mathcal{P}_{\kappa,\ell}|<|\mathrm{i}\mathcal{Q}_{\kappa,\ell}|$ on the circle, so $\mathcal{D}_{\kappa,\ell}$ also has exactly one zero there. The constant in $\rho_\ell$ is $C=4B_2/c_1$.
\end{proof}
\begin{proof}[Proof of Theorem \ref{maintheorem2}]
Let $\kappa=\eta_1/\eta_0$ and let $z_\ell$ ($\ell\ge\ell_0$) be the zeros given by Proposition \ref{wgm:rouche}. By \eqref{wgm:DD}, $\mathcal{D}^1_\ell(z_\ell)=0$. Set $\lambda_\ell=z_\ell/(r_0\sqrt{\epsilon_0\mu_0})$.
Put $\Lambda_\ell=\Lambda(x_\ell/\nu)$. Lemma~\ref{wgm:debye}(i) and the Cauchy bounds in the proof of Proposition~\ref{wgm:rouche} give, on $|z-x_\ell|\le\rho_\ell$,
\begin{align*}
|\vartheta_\ell(z)|&\ge(c-C\rho_\ell\nu^{1/2})\mathrm{e}^{\nu\Lambda_\ell},
&|\psi_\ell(z)|&\le C\nu^{1/2}\mathrm{e}^{-\nu\Lambda_\ell}.
\end{align*}
For large $\ell$, the first bound exceeds the second, so $\xi_\ell(z)=\psi_\ell(z)+\mathrm{i}\vartheta_\ell(z)\ne0$ throughout this disc.
\begin{enumerate}
\item By Proposition \ref{noupper} and its proof, $\mathcal{D}^1_\ell(z_\ell)=0$ gives the non-trivial purely outgoing solutions $u_{\ell m}=t\,M^0_{\ell m}(k_1r)$ ($r<r_0$), $u_{\ell m}=M^-_{\ell m}(k_0r)$ ($r>r_0$), with $t=h^{(1)}_\ell(x_0)/j_\ell(x_1)$ ($j_\ell(x_1)\neq0$, since otherwise $\mathcal{D}^1_\ell(z_\ell)=\eta_0\xi_\ell(z_\ell)\psi_\ell'(m_0z_\ell)\neq0$). Hence $\lambda_\ell$ is a resonance of $P_E$, with the $2\ell+1$ resonance states $u_{\ell m}$, $|m|\le\ell$; it is a genuine pole of $A^1_\ell$, because $N^1_\ell$ and $D^1_\ell$ have no common zero (if both vanish, the Wronskian $\psi_\ell\vartheta_\ell'-\psi_\ell'\vartheta_\ell=1$ at $x_0$ forces $\psi_\ell(x_1)=\psi_\ell'(x_1)=0$).
\item Proposition~\ref{nohalfplane} gives $\Im\lambda_\ell<0$.
\item Since $|z_\ell-x_\ell|<\rho_\ell\le1$ and $x_\ell\in(t_1\nu,t_2\nu)$, we have $t_1\nu-1<\Re z_\ell<t_2\nu+1$, so $\Re\lambda_\ell\to\infty$.
\item As $\Lambda$ is decreasing, $|\Im z_\ell|<\rho_\ell\le C\nu^{\frac32}\mathrm{e}^{-2\nu\Lambda(t_2)}\le C'\mathrm{e}^{-\nu\Lambda(t_2)}$. Using $\nu>(\Re z_\ell-1)/t_2$,

\begin{equation*}
|\Im\lambda_\ell|\le C\mathrm{e}^{-c\Re\lambda_\ell},\qquad c=\frac{r_0\sqrt{\epsilon_0\mu_0}\,\Lambda(t_2)}{t_2}>0 .
\end{equation*}

The constants depend only on $m_0,\eta_0,\eta_1,r_0$ and $\epsilon_0\mu_0$.
\item Given $C_1,C_2>0$, (iii) and (iv) show that $\lambda_\ell\in\{|\Re\lambda|\ge C_1,\ \Im\lambda\ge-C_2\}$ for all large $\ell$, so no resonance free strip exists.
\end{enumerate}
\end{proof}

\section*{Acknowledgements}
The author A.W. wishes to thank Yan-Long Fang and Nelia Charalambous for useful discussions about Mie scattering that led to this paper. She also wishes to thank Andras Vasy for suggesting the topic of whispering gallery modes.

\section*{Declarations}

\begin{itemize}
\item \textbf{Funding.} None
\item \textbf{Competing interests.} None
\item \textbf{Data availability.} Data sharing is not applicable to this article, as no datasets were generated or analysed during the current study.
\item \textbf{Ethics approval and consent to participate.} Not applicable.
\item \textbf{Consent for publication.} Not applicable.
\item \textbf{Materials availability.} Not applicable.
\item \textbf{Code availability.} Not applicable.
\item \textbf{Author contribution.} Not applicable (single author).
\item \textbf{Use of large language models.} All ideas here are the original work of the author. Claude Opus 5.5 AI was used for proofreading/language organization and for assistance to develop the difficult Bessel function asymptotics specifically in Regions 2 and 3 of Theorem \ref{maintheorem} and Theorem \ref{maintheorem2}, after the various regimes were fixed first by the author. 
\end{itemize}

\begin{appendices}

\section{Outgoing asymptotics and uniqueness}
The outgoing condition in Definition~\ref{resonanceP} depends only on the field near infinity. For real $\lambda\ne0$, an outgoing field satisfying the free equation outside a ball has the expansion
\begin{align}\label{OR}
f=\frac{\mathrm{e}^{\mathrm{i} k_0 r}}{r}\left(\Phi+\mathcal{O}\left(\frac{1}{r}\right)\right) \quad r\rightarrow \infty
\end{align}
with $\Phi\in C^{\infty}(\mathbb{S}^{2};\mathbb C^3)$ the restriction of an entire function on $\mathbb C^3$ to the sphere. The expansion can be differentiated in $r$.

\begin{theorem}[Maxwell transmission uniqueness]\label{Rellich}
Let $\Omega$ be any bounded Lipschitz open set with the positive coefficients \eqref{eq:lipschitz-coefficients}. Let $\lambda\in\mathbb R\setminus\{0\}$. An outgoing $u\in H_{\loc}(\curl)$ satisfying $(P_E-\lambda^2)u=0$ and the transmission conditions \eqref{transmission} vanishes. The same holds for $P_H$.
\end{theorem}
\begin{proof}[Sketch adapted from \cite{DZbook}]
The equation gives $\div(\epsilon u)=0$, so the outgoing coefficient $\Phi$ in \eqref{OR} is tangential. Green's identity, with cancellation of the interface terms, yields
\begin{equation*}
0=\lim_{R\to\infty}\Im\int_{\partial B_R}
(\hat r\times\mu_0^{-1}\curl u)\cdot\overline u\,\mathrm{d}S
=-\frac{k_0}{\mu_0}\|\Phi\|_{L^2(\mathbb S^2)}^2.
\end{equation*}
Thus $\Phi=0$, and the scalar Rellich lemma applied componentwise gives $u=0$ in the unbounded component of $\mathbb R^3\setminus\overline\Omega$. Set $H=(\mathrm{i}\lambda\mu)^{-1}\curl u$. At an interface neighborhood adjoining a region where both fields vanish, their tangential traces on the other side are zero, and the fields can be extended by zero. These extensions solve the constant-coefficient Maxwell system for that side in a full neighborhood of the patch. Their components solve the Helmholtz equation distributionally and are analytic there; vanishing on one side forces vanishing across the neighborhood. Interior unique continuation and repetition across the finitely many interface components give $u=H=0$ everywhere, including any bounded complementary components. Interchanging $\epsilon,\mu$ for $P_H$, the result also holds for the magentic operator.  Conjugation also gives incoming uniqueness.
\end{proof}

\section{Bessel and Hankel estimates}\label{app:bessel}
We use $\nu$, $\psi_\ell$ and $\xi_\ell$ from \eqref{eq:riccati-definitions}; primes denote differentiation with respect to the full argument.

\subsection{Basic identities}\label{AA}
The Riccati equation is \eqref{wgm:ode}. For $f_\ell=j_\ell,h_\ell^{(1)}$ or $h_\ell^{(2)}$, the derivative recurrence is
\begin{align}\label{recurrence}
f_\ell'(z)=f_{\ell-1}(z)-\frac{\ell+1}{z}f_\ell(z).
\end{align}
The spherical Bessel Wronskian is
\begin{align}\label{wronskian}
j_{\ell}(z)y_{\ell}'(z)-y_{\ell}(z)j_{\ell}'(z)=z^{-2}.
\end{align}

\subsection{Subturning estimates}\label{app:subturning}
The real Debye estimates away from the turning point are stated in Lemma~\ref{wgm:debye}. We record here the uniform logarithmic-derivative estimates needed up to the edge of the turning window.

For Region~1, fix $a,C_0>0$ and assume

\begin{equation*}
a\leq\Re w\leq\nu-M\nu^{1/3},\quad |\Im w|\leq C_0,
\quad M\geq M_0,\quad \nu\geq\nu_0(M,a,C_0).
\end{equation*}

Set $\beta(w)=\sqrt{\nu^2-w^2}/w$, positive for $0<w<\nu$, and $\rho_\nu(w)=\nu^2/|\nu^2-w^2|^{3/2}$. The uniform Bessel expansions and their derivatives, \cite{NIST}, (10.20.4), (10.20.5), (10.20.7), (10.20.8) and \S10.20(ii), combined with the Airy expansions (9.7.5)--(9.7.8), give
\begin{align*}
\frac{J_\nu'(w)}{J_\nu(w)}&=\beta(w)(1+O(\rho_\nu(w))),\\
\frac{Y_\nu'(w)}{Y_\nu(w)}&=-\beta(w)(1+O(\rho_\nu(w))),\\
\frac{(H_\nu^{(1)})'(w)}{H_\nu^{(1)}(w)}&=-\beta(w)(1+O(\rho_\nu(w))).
\end{align*}
The constants depend only on $a,C_0$. Since $|\nu^2-w^2|\geq cM\nu^{4/3}$, the relative errors are $O(M^{-3/2})$, uniformly at both complex arguments $w=z,m_0z$ used in Region~1.

\subsection{Two-term Airy expansions in the turning windows}\label{airy-two-term}
For the full complex argument $w=\nu+\tau\nu^{1/3}$, put
\begin{gather*}
s=-2^{1/3}\tau,\qquad p=2^{1/3}\mathrm{e}^{-\pi\mathrm{i}/3},\qquad \sigma=p\tau,\\
A=\operatorname{Ai}(s),\qquad \widehat A=\operatorname{Ai}(\sigma).
\end{gather*}
Here $A'$ and $\widehat A'$ denote derivatives of the Airy function with respect to its argument. We use $w=z$ in the exterior turning window and $w=m_0z$ in the interior turning window.

\begin{proposition}[Uniform two-term Airy expansion]\label{prop:two-term-airy}
Fix $M>0$ and $0<\delta<\pi$. There are constants $C_{M,\delta}>0$ and $\nu_0=\nu_0(M,\delta)$ such that, for every positive half-integer $\nu\geq\nu_0$ and every $|\tau|\leq M$, with $|\operatorname{ph}(w/\nu)|\leq\pi-\delta$, the following formulas hold:
\begin{equation}\label{eq:airy-psi}
\psi_\ell(w) ={} \sqrt\pi\,2^{-1/6}\nu^{1/6} \Big\{A+\frac{1}{\nu^{2/3}} \Big(\frac{3\tau}{10}A +\frac{3\,2^{1/3}\tau^2}{10}A'\Big)\Big\} +r_\psi,
\end{equation}
\begin{equation}\label{eq:airy-psiprime}
\psi_\ell'(w) ={} \sqrt\pi\,2^{-1/6}\nu^{-1/6} \Big\{-2^{1/3}A' +\frac{1}{\nu^{2/3}} \Big[\Big(\frac{3}{10}+\frac{3\tau^3}{5}\Big)A +\frac{3\,2^{1/3}\tau}{10}A'\Big]\Big\} +r_{\psi'},
\end{equation}
\begin{equation}\label{eq:airy-xi}
\xi_\ell(w) ={} \sqrt\pi\,2^{5/6}\mathrm{e}^{-\pi\mathrm{i}/3}\nu^{1/6} \Big\{\widehat A+\frac{1}{\nu^{2/3}} \Big(\frac{3\tau}{10}\widehat A +\frac{3\,2^{1/3}\mathrm{e}^{2\pi\mathrm{i}/3}\tau^2}{10}\widehat A' \Big)\Big\} +r_\xi,
\end{equation}
\begin{equation}\label{eq:airy-xiprime}
\xi_\ell'(w) ={} \sqrt\pi\,2^{5/6}\mathrm{e}^{-\pi\mathrm{i}/3}\nu^{-1/6} \Big\{p\widehat A' +\frac{1}{\nu^{2/3}} \Big[\Big(\frac{3}{10}+\frac{3\tau^3}{5}\Big)\widehat A -\frac{3p\tau}{10}\widehat A'\Big]\Big\} +r_{\xi'},
\end{equation}
The remainders satisfy the uniform absolute bounds
\begin{align}
 |r_\psi|+|r_\xi|&\leq C_{M,\delta}\nu^{-7/6},
 &
 |r_{\psi'}|+|r_{\xi'}|&\leq C_{M,\delta}\nu^{-3/2}.
\label{eq:airy-riccati-errors}
\end{align}
\end{proposition}
\begin{proof}
Apply the transition expansions (10.19.8), (10.19.9), (10.19.12) and (10.19.13) of \cite{NIST} to \eqref{eq:riccati-definitions}, and expand $(w/\nu)^{1/2}=1+\frac12\tau\nu^{-2/3}+O_M(\nu^{-4/3})$. This gives the displayed coefficients.

The complex uniform expansions in \S10.20(ii) of \cite{NIST}, with the remainder estimates on pp.~419--425 of \cite{olverbook}, give absolute errors $O_{M,\delta}(\nu^{-5/3})$ for $J_\nu,H_\nu^{(1)}$ and $O_{M,\delta}(\nu^{-2})$ for their derivatives on $|\tau|\leq M$. Multiplying by $\sqrt{\pi w/2}$ and including the derivative of this prefactor yields \eqref{eq:airy-riccati-errors}.
\end{proof}
These additive estimates remain valid at zeros of $A$ or $A'$; no division by $\psi_\ell(m_0z)$ is needed in the interior turning window.

\subsection{Exterior Airy expansion}\label{largeturning}

For $q>0$ we define $\zeta=\zeta(q)$ by
\begin{align}
\left(\frac{\mathrm{d}\zeta}{\mathrm{d}q}\right)^2=\frac{1-q^2}{\zeta q^2}.
\end{align}
This is piecewise defined for real $q$ as
\begin{align}
\frac{2}{3}\zeta^{\frac{3}{2}}=\log\left(\frac{1+\sqrt{1-q^2}}{q}\right)-\sqrt{1-q^2} \quad 0<q\leq 1
\end{align}
and
\begin{align}\label{eq:zeta-oscillatory}
\frac{2}{3}(-\zeta)^{\frac{3}{2}}=\sqrt{q^2-1}-\mathrm{arcsec} q \quad 1\leq q<\infty.
\end{align}
These are equations (10.20.1)--(10.20.3) of \cite{NIST}.
For $q=z/\nu$ in Region~3, the uniform Hankel expansions \cite[(10.20.6), (10.20.9), \S10.20(ii)]{NIST} give \eqref{eq:exterior-logderivative}, with the uniform $O(\nu^{-1})$ remainder stated there. The range is
\begin{equation*}
1+M\nu^{-2/3}\leq\Re q<K\nu/m_0,\qquad |\Im q|\leq C_2/\nu,
\end{equation*}
with branches continued from real $q>1$. Only this logarithmic-derivative estimate is needed. The Airy phase convention follows from
$\operatorname{Ai}(s)-\mathrm{i}\operatorname{Bi}(s)=2\mathrm{e}^{-\pi\mathrm{i}/3}\operatorname{Ai}(\mathrm{e}^{2\pi\mathrm{i}/3}s)$.

\subsection{Large-argument estimates}
Set
\begin{equation}\label{eq:large-phase}
\omega(w)=w-\frac{\pi\nu}{2}-\frac{\pi}{4}.
\end{equation}
For $|w|\geq K\nu^2$, $K\geq1$, and $|\Im w|\leq C_0$, we have
\begin{align}
\psi_\ell(w)&=\cos\omega(w)+O\left(\frac{\nu^2}{|w|}\right),\label{eq:largepsi}\\
\psi_\ell'(w)&=-\sin\omega(w)+O\left(\frac{\nu^2}{|w|}\right),\label{eq:largepsip}\\
\xi_\ell(w)&=\mathrm{e}^{\mathrm{i}\omega(w)}\left(1+O\left(\frac{\nu^2}{|w|}\right)\right),\label{eq:largexi}\\
\xi_\ell'(w)&=\mathrm{e}^{\mathrm{i}\omega(w)}\left(\mathrm{i}+O\left(\frac{\nu^2}{|w|}\right)\right).\label{eq:largexip}
\end{align}
The constants are uniform in $\ell\geq1$ and $K\geq1$ and depend only on $C_0$. The estimates for $\xi_\ell,\xi_\ell'$ also hold without the restriction on $\Im w$; for fixed $\ell$ their remainders are $O_\ell(|w|^{-1})$ as $|w|\to\infty$.

\begin{proof}
Write $\xi_\ell(w)=\mathrm{e}^{\mathrm{i}\omega(w)}(1+r_\ell(w))$. In the finite spherical Hankel formula \cite[(10.49.1), (10.49.6)]{NIST}, the coefficient of $w^{-k}$ has modulus at most $[\ell(\ell+1)/2]^k/k!$. Thus, for $w\ne0$,
\begin{align*}
|r_\ell(w)|&\leq\exp\left(\frac{\ell(\ell+1)}{2|w|}\right)-1,\\
|r_\ell'(w)|&\leq\frac{\ell(\ell+1)}{2|w|^2}\exp\left(\frac{\ell(\ell+1)}{2|w|}\right).
\end{align*}
Since $\ell(\ell+1)<\nu^2$, these imply \eqref{eq:largexi}--\eqref{eq:largexip}. The identity $\psi_\ell(w)=\frac12(\xi_\ell(w)+\overline{\xi_\ell(\overline w)})$ and its derivative give \eqref{eq:largepsi}--\eqref{eq:largepsip}, using $|\mathrm{e}^{\pm\mathrm{i}\omega(w)}|\leq\mathrm{e}^{C_0}$.
\end{proof}

\end{appendices}

\textbf{AI usage}. Claude Opus 5.5 was used in the development of the estimates in Regions 2 and 3 in both main theorems. 

\begingroup
\fontsize{10}{12}\selectfont

\endgroup

\end{document}